\documentclass[english]{article}
\usepackage[T1]{fontenc}
\usepackage[latin9]{inputenc}
\usepackage{babel}
\usepackage{amsfonts,color}

\newtheorem{dfn}{Definition}[section]
\newtheorem{theorem}[dfn]{Theorem}

\newtheorem{remark}[dfn]{Remark}

\newtheorem{lem}[dfn]{Lemma}

\newtheorem{cor}[dfn]{Corollary}
\newtheorem{fact}[dfn]{Fact}
\usepackage{amssymb} 
\usepackage{amsmath}

 \global\long\def\sbr#1{\left[ #1\right] }
 \global\long\def\cbr#1{\left\{  #1\right\}  }
 
 \global\long\def\rbr#1{\left(#1\right)}
\global\long\def\ev#1{\mathbb{E}{#1}}
 \global\long\def\E{\mathbb{E}}
 \global\long\def\P{\mathbb{P}}
 \global\long\def\R{\mathbb{R}}
  \global\long\def\N{\mathbb{N}}
  
  \global\long\def\Z{\mathbb{Z}}

 \global\long\def\conv{\rightarrow}

 \global\long\def\dd#1{\textnormal{d}#1}

 \global\long\def\TTV#1#2#3{\text{TV}^{#3}\!\rbr{#1,#2}}
 
 \global\long\def\UTV#1#2#3{\text{UTV}^{#3}\!\rbr{#1,#2}}
 \global\long\def\DTV#1#2#3{\text{DTV}^{#3}\!\rbr{#1,#2}}
 \global\long\def\TTVemph#1#2#3{\emph{TV}^{#3}\!\rbr{#1,#2}}
 \global\long\def\UTVemph#1#2#3{\emph{UTV}^{#3}\!\rbr{#1,#2}}
 \global\long\def\DTVemph#1#2#3{\emph{DTV}^{#3}\!\rbr{#1,#2}}
 \global\long\def\cad{c\`{a}dl\`{a}g }

 \global\long\def\ra{\rightarrow}
 
 \global\long\def\ns{\infty}
 \global\long\def\8{\infty}

 \global\long\def\ns{\infty}
\global\long\def\Ucross#1#2#3{\text{u}^{#1}\!\rbr{#2,#3}}
\global\long\def\Dcross#1#2#3{\text{d}^{#1}\!\rbr{#2,#3}}
\global\long\def\cross#1#2#3{\text{n}^{#1}\!\rbr{#2,#3}}

\global\long\def\Ucrossemph#1#2#3{\emph{u}^{#1}\!\rbr{#2,#3}}
 \global\long\def\Dcrossemph#1#2#3{\emph{d}^{#1}\!\rbr{#2,#3}}
 \global\long\def\crossemph#1#2#3{\emph{n}^{#1}\!\rbr{#2,#3}}
 \global\long\def\UTVemph#1#2#3{\emph{UTV}^{#3}\!\rbr{#1,#2}}
\global\long\def\DTVemph#1#2#3{\emph{DTV}^{#3}\!\rbr{#1,#2}}

\title{Local times of deterministic paths and self-similar processes with stationary increments as normalized numbers of interval crossings}

\author{Witold Bednorz, Purba Das and Rafa{\l } {\L }ochowski}

\begin{document}

\maketitle

\begin{abstract}
We prove a general result on a relationship between a limit of
normalized numbers of interval crossings by a c\`adl\`ag path and an occupation
measure associated with this path.

Using this result, we define local times of fractional Brownian motions (classically defined as densities of relevant occupation measure) as weak limits of properly normalized numbers of interval crossings.

We also discuss a similar result for c\`adl\`ag semimartingales,
in particular for alpha-stable processes, and Rosenblatt processes, and provide natural examples
of deterministic paths which possess quadratic or higher-order variation but no local
times.
\end{abstract}

\section{Introduction}
In this article, we deal with local times of real c\`adl\`ag paths which may be defined as weak limits of normalized numbers of interval crossings by these paths. Instances of real c\`adl\`ag paths possessing such local times will be trajectories of self-similar stochastic processes with stationary increments, like fractional Brownian motions or alpha-stable processes.

The idea to define the local time of a standard Brownian motion as a limit of (properly normalized) numbers of (properly defined) interval crossings by trajectories of this process comes already from Paul L\'evy \cite{Levy:1940}. This definition was later generalized to the case of continuous semimartingales \cite{ElKaroui:1978}, \cite[Chapt. VI]{RevuzYor:2005} and Markov processes \cite{FristedTaylor:1983} (to be more precise, in \cite{FristedTaylor:1983} the Authors deal rather with excursions from a given level than with the interval crossings, but both approaches coincide for processes with a.s. continuous paths). Since the number of interval crossings may be calculated for each path separately, such results provide a pathwise construction of local time. Results of this type for a standard Brownian motion  $B$ (and continuous semimartingales) were obtained first by Chacon et al. \cite{Chacon:1981}, built on ideas from \cite{Perkins:1981}.  In \cite{Chacon:1981} the authors prove the existence of a measurable set $\Omega_0$ such that ${\Omega_0}$ has probability $1$ and for each $\omega \in \Omega_0$ there exists the mentioned limit of normalized numbers of interval crossings by the trajectory $B(\omega)$. The excellent Master thesis of Marc Lemieux \cite{Lemieux:1983} extended this result to the case of c\`adl\`ag semimartingales with locally summable jumps. The case of general c\`adl\`ag semimartingales was treated only recently in \cite{LochOblPS:2021}.

Interval crossings by deterministic or random paths received attention in many recent publications on pathwise stochastic calculus or model-free finance, see for example \cite{PerkowskiProemel_local:2015, Obloj_local:2015, ContPerkowski:2018, LochOblPS:2021, Kim:2022}. 

Given most of the above-mentioned results are obtained for semimartingales and Markov processes, a natural question arises whether similar results hold for processes possessing local times (in the sense defined below), which do not belong to these classes. Examples of such processes are fractional Brownian motions. Fractional Brownian motion (fBm in short) possesses local time, but its characterization as a limit of normalized numbers of interval crossings by its trajectories has not been much investigated so far. This seems surprising considering the age of L\'evy's result. A recent paper on this approach for fBms with the Hurst index less than $1/2$ is \cite{Toyomu:2023}. In this note, we consider such results for self-similar processes with stationary increments. In particular, we obtain the existence of the limit of normalized numbers of interval crossings for fBms with a Hurst index from the whole interval $(0,1)$. However, the convergence we obtain is the weak convergence of measures, thus much weaker than the almost sure uniform convergence of processes obtained in \cite{Toyomu:2023}. 

Another example of self-similar processes with stationary increments are $\alpha$-stable processes. They are also special cases of c\`adl\`ag semimartingales considered in \cite{LochOblPS:2021}, however, \cite{LochOblPS:2021} considers the semimartingale local time (defined as the density of the occupation measure along the 'business clock' -- the continuous part of the quadratic variation); here we consider the local time defined as the density of the occupation measure along the `natural clock' -- the Lebesgue measure on $[0, +\ns)$. 

\textbf{Preview}. In the next section, we state the main results obtained in this paper, the proofs of these results are presented in subsequent sections. In Section \ref{sectwo}, we present essential definitions, notations, and prove auxiliary results (Lemma \ref{Ban_Viatli_extended} and Theorem \ref{thm:meta1}), which may be of independent interest. Examples of application of Theorem \ref{thm:meta1} to deterministic paths are presented in Subsection \ref{subsect:examples}. In Section \ref{Secthree} we apply Theorem \ref{thm:meta1} to paths of self-similar processes with stationary increments. In Section \ref{secfour} we list examples of self-similar processes with stationary increments for which the results of Section \ref{Secthree} may be applied. They include, apart from the mentioned fBms and $\alpha$-stable processes, also the Rosenblatt processes. We also present several examples where assumptions of Theorem \ref{thm:meta1} hold, but no local time exists. 

\section{Main results} \label{sectwo0}
To state the main results, we need to start with several definitions.

\begin{dfn} \label{tvc_def} The truncated variation of $x: [0, +\ns) \ra \R$ with the truncation parameter $c \ge 0$ on the time interval $[s,t]$, $0 \le s < t < +\ns$, is defined as:
\begin{equation} \label{tvc_def-eq}
\TTVemph x{\left[s,t\right]}{c}:=\sup_{\pi \in \Pi(s,t)}\sum_{[u,v] \in \pi} \rbr{\left| x_{v}-x_{u} \right| -c}_{+},
\end{equation}
where $(\cdot)_+ = \max\rbr{\cdot, 0}$ and the supremum is taken over all finite partitions $\pi$ of the interval $[s,t]$, that is finite sets of no overlapping (with disjoint interiors) subintervals $[u,v]$ of $[s,t]$ such that $\bigcup_{[u,v] \in \pi}[u,v] = [s,t]$. The family of all such partitions is denoted by $\Pi(s,t)$.

Similarly, the upward and downward truncated variations of $x$ are defined respectively as 
\begin{equation}\label{eq:uptrunc variation}
  \UTVemph {x}{[s,t]}{c} := \sup_{\pi \in \Pi(s,t)}\sum_{[u,v] \in \pi} \rbr{ x_{v}-x_{u} -c}_+
\end{equation}
and
\begin{equation}\label{eq:downtrunc variation}
   \DTVemph{x}{[s,t]}{c} := \sup_{\pi \in \Pi(s,t)}\sum_{[u,v] \in \pi} \rbr{ x_{u}-x_{v} -c}_+.
\end{equation}   
\end{dfn}

The truncated variation is finite for any c\`adl\`ag or even regulated (i.e. possessing right- and left- limits) path $x$ whenever $c>0$. If $c=0$ then the truncated variation coincides with the \emph{total variation of $x$} while the upward and downward truncated variations coincide with the \emph{upward total variation} and \emph{downward total variation}. They are denoted as $ \TTV {x}{[s,t]}{}$, $ \UTV {x}{[s,t]}{}$ and $ \DTV {x}{[s,t]}{}$ respectively.

Let $\cross{y ,c}{x}{[0, t]}$,  $\Ucross{y ,c}{x}{[0, t]}$ and $\Dcross{y ,c}{x}{[0, t]}$ be the numbers of \emph{crossings, upcrossings and downcrossings} by $x$ the interval $[y-c/2, y+c/2]$ on the time interval $[s,t]$ respectively. For precise definitions, we refer to Definition \ref{defd}. 
The following result, which is of independent interest, is one of the main ingredients we will use in establishing the weak convergence of normalized numbers of interval crossings to local times. 
{As far as we know, no similar result has been proven in literature so far.}
\begin{theorem} \label{thm:meta1} 
Let $x,\zeta:[0 ,+\infty) \rightarrow \mathbb{R}$ be c\`adl\`ag, and $\varphi:(0 ,+\infty) \rightarrow[0 ,+\infty)$ be  
such that $\lim _{c \downarrow 0} \varphi(c)=0$. Assume that 
\begin{align}
\label{conv_tv0}
\forall t \geq 0, \quad \varphi(c) \mathrm{TV}^{c}(x,[0, t]) \rightarrow \zeta_{t} \quad \text{ as } c \downarrow 0,
\end{align}
or equivalently that any of the following equivalent conditions holds
\begin{align}
\forall t \geq 0, \quad 2\varphi(c) \mathrm{UTV}^{c}(x,[0, t]) \rightarrow \zeta_{t} \quad \text{ as } c \downarrow 0,  \label{conv_utv} \\
\forall t \geq 0, \quad 2\varphi(c) \mathrm{DTV}^{c}(x,[0, t]) \rightarrow \zeta_{t} \quad \text{ as } c \downarrow 0. \label{conv_dtv}
\end{align}
Then $\zeta$ is non-negative, non-decreasing and for any continuous $g: \mathbb{R} \rightarrow \mathbb{R}$ we have the following pointwise
convergence of Lebesgue-Stieltjes integrals
\begin{align}
\label{weak_local_tv}
\forall t \geq 0,  \quad \varphi(c) \int_{\mathbb{R}} \mathrm{n}^{y, c}(x,[0, t]) g(y) \mathrm{d} y \rightarrow \int_{(0,t]} g\left(x_{s-}\right) \mathrm{d} \zeta_{s} \quad \text{as }  c \downarrow 0.
\end{align}
The convergence in \eqref{weak_local_tv}  also holds when $\mathrm{n}^{z, c}(x,[0, t])$ is replaced by $2\mathrm{u}^{z, c}(x,[0, t])$ or by $2\mathrm{d}^{z, c}(x,[0, t])$, and if $\zeta$ is continuous such convergence is uniform on compacts. 
\end{theorem}
{For the proof of Theorem \ref{thm:meta1} we refer to Sect. \ref{aux_res}.}
\begin{remark} \label{normalisation}
If $\TTVemph x{[0,1]}{} = +\ns$ then defining the normalization function by:
\[
\varphi(c) := \frac{1}{1+\TTVemph x{[0,1]}{c}}
\] 
we obtain a non-decreasing, non-negative function such that $\lim_{c\ra 0+}\varphi\left(c\right)=0$, $\lim_{c\ra +\ns}\varphi\left(c\right)= 1$ and 
\[
\zeta_1 = \lim_{c\ra 0+}\varphi\left(c\right) \TTVemph x{[0,1]}{c}= 1.
\]
\end{remark} 

Relation \eqref{weak_local_tv} may be seen as a `weak' form of occupation times formula (see the next section).  The main difference between both formulae is that $g$ in \eqref{weak_local_tv} is assumed to be \emph{continuous} and that we do not know whether quantities $\varphi(c)  \cross{y,c}{x}{[0, t]}$ tend in some sense to any limit; we only know (and this may be an equivalent statement of the thesis) that the measures $\varphi(c) n^{y,c}\left(x,\left[0, t \right]\right) \dd y$ tend weakly to the occupation measure
\[
\mu_t(\Gamma) = \int_{(0,t]} {\mathbf 1}_{\Gamma}(x_{s-}) \zeta(\dd s).
\] However, if \eqref{conv_tv0}, \eqref{conv_utv} or \eqref{conv_dtv} holds, $\zeta$ is continuous and $x$ possesses the local time $L$ relative to the measures $\dd y$ (the Lebesgue measure on $\R$) and $\zeta$, then  the following relations 
\[
\int_{(0, t]}g(x_{s-})\zeta(\dd s) = \int_{(0, t]}g(x_{s})\zeta(\dd s) = \int_{\R} g(y) L_t^y \dd y
\]
hold for any continuous $g: \R \ra \R$, which together with \eqref{weak_local_tv} implies  that for any finite $t>0$, $L_t^y \dd y$ is the weak limit of the measures $\varphi(c) \cross{y,c}{x}{[0, t]} \dd y$. In one of the examples (see Sect.  \ref{subsect46}), we show that this convergence can not be replaced by a stronger mode of convergence -- the weak convergence of $\varphi(c) \cross{y,c}{x}{[0, t]}$ to $L_t^y$ in $\mathbb{L}^1 (\R, \dd y)$.

For some processes it may be proven that for each $t \ge 0$ the convergence
\begin{equation} \label{convas}
c^{1/\beta-1}\TTV X{[0, t]}c\ra C \cdot t \text{ as } c \ra 0+
\end{equation} 
holds with probability $1$, like for example for fBms with arbitrary Hurst parameter $H \in(0,1)$, see Sect. \ref{secfour}. Then, as a direct consequence of Theorem \ref{thm:meta1} applied to each trajectory of $X$ separately, we have the following result.
\begin{cor}\label{cor:almost_sure}
Assume that for a real process $X_t$, $t \in [0, +\ns)$, the convergence in \eqref{convas} holds with probability $1$ and that $X$ possesses a local time $L$  (in the sense of Definition \ref{localt_process}) relative to the Lebesgue measure $\dd y$ on $\R$ and the constant mapping $\Xi(\omega) \equiv \dd u$ ($\dd u$  is the Lebesgue measure on $[0, +\ns)$).  Then there exist a measurable set $\Omega_X \subset \Omega$ such that $\P \rbr{\Omega_X} = 1$ and for each $\omega \in \Omega_X$ and each $t\in[0, +\ns)$
\[
c^{1/\beta-1} \crossemph{y ,c}{X(\omega)}{[0, t]} \dd y  \ra^{weakly} C \cdot L_t(\omega) \dd y \text{ as } c \ra 0+,
\]
where $C$ is the same as in  \eqref{convas}  and `$\ra^{weakly}$' denotes the weak convergence of measures.
\end{cor}

Existence of the deterministic limit \eqref{conv_tv0} of the form $\zeta(t) = C \cdot t$  (for a proper normalization $\varphi(\cdot)$ and a deterministic constant $C$) may be established for a class of self-similar processes with stationary increments with the help of subadditive ergodic theorem and the convergence  holds in probability. 
\begin{theorem}\label{thm:ssProcesses} Let $\rbr{\Omega, {\cal F}, \P}$ be a probability triple and $X_t$, $t \in [0,+\ns)$, be a self-similar process with index $\beta >0$ on $\rbr{\Omega, {\cal F}, \P}$, which has stationary increments and \cad trajectories. We assume additionally that for some $c_0>0$ it fulfils 
	\begin{equation}
		\ev{}\TTVemph X{[0;1]}{c_0}<+\infty \label{eq:finiteExpectation}
\end{equation}
and that for any $k \in \N$ the stationary sequence
$$\rbr{\TTVemph X{[0,k]}1, {\TTVemph X{[k,2k]}1}, {\TTVemph X{[2k,3k]}1}, \ldots }$$
is ergodic.
Then there exists $C \in (0, +\ns)$ such that  for each $t\in[0, +\ns)$ the following convergence in probability holds
\begin{equation} \label{convergodic}
c^{1/\beta-1}\TTVemph X{[0, t]}c\ra^{\P} C \cdot t \text{ as } c \ra 0+.
\end{equation}
\end{theorem} 
{For the proof of Theorem \ref{thm:ssProcesses} and next  Corollary \ref{cor:ltssProcesses} we refer to Sect. \ref{Secthree}.}
\begin{remark} \label{rem:constant} The main drawback of Theorem \ref{thm:ssProcesses} is
that it does not explicitly identify  the constant $C$, however, from its proof it follows that 
\[
C = \lim_{n \ra +\ns} \frac{\E \TTVemph X{[0,n]}1}{n}.
\]
From the fact that $\rbr{\TTVemph X{[m,n]}1+1}$ is subadditive ($\TTVemph X{[0,n]}1 +1 \le \TTVemph X{[0,m]}1+1+\TTVemph X{[m,n]}1+1$), see \eqref{superaddd}, we get even tighter estimate: for any $n \in \N$
\[
\frac{\E \TTVemph X{[0,n]}1}{n} \le C \le \frac{\E \TTVemph X{[0,n]}1+1}{n}.
\]
\end{remark}
\begin{cor}\label{cor:ltssProcesses}
Assume that a real process $X_t$, $t \in [0, +\ns)$, satisfies assumptions of Theorem  \ref{thm:ssProcesses} and possesses a local time $L$  relative to the Lebesgue measure $\dd y$ on $\R$ and the constant mapping $\Xi(\omega) \equiv \dd u$.  Then there exists a sequence $\rbr{c_n}$, $n\in \N$, such that $c_n \ra 0+$ as $n \ra +\ns$ and a measurable set $\Omega_X \subset \Omega$ such that $\P \rbr{\Omega_X} = 1$ and for each $\omega \in \Omega_X$ and  each $t\in[0, +\ns)$
\[
c_n^{1/\beta-1} \crossemph{y ,c_n}{X(\omega)}{[0, t]} \dd y \ra^{weakly} C \cdot L_t(\omega) \dd y \text{ as } n \ra +\ns,
\]
where $C$ is the same as in Remark  \ref{rem:constant}. 
\end{cor}
\noindent Examples of application of Corollaries  \ref{cor:almost_sure} and \ref{cor:ltssProcesses}  to specific processes are presented in Sect.  \ref{secfour}.

\section{Notation, definitions and auxiliary results} \label{sectwo}
\textit{Notations:} By $\Z$ we denote the set of all integers, by $\N$ we denote the set of all positive integers and by $\N_0$ we denote the set $\N \cup \cbr{0}$. By $C([0,+\ns);\R)$ we denote the space of all continuous functions $x\colon [0,+\ns) \ra \R$.
By $D([0,+\ns);\R)$ we denote the space of all c{\`a}dl{\`a}g (RCLL) functions $x\colon [0,+\ns) \ra \R$, that is $x \in D([0,+\ns);\R)$ if it is right-continuous at each $t\in [0,+\ns)$ and possesses left-limits at each $t\in (0,+\ns)$. 

For $x\in D([0,+\ns);\R)$ we set $x_{t-}:=\lim_{s<t, s\to t}x_t$ for $t\in (0,+\ns)$, $x_{0-}:=x_0$ and $\Delta x_s := x_s-x_{s-}$ for $s\in [0,+\ns)$.

$V^0([0,+\ns);\R)$ denotes the subset of $D([0,+\ns);\R)$ of piecewise monotonic \cad paths, that is functions $x\in D([0,+\ns);\R)$ for which for any $T>0$ there exist \emph{finite} number of intervals $I_i$, $i =1,2,\ldots, N$, $N \in \N$, such that $\bigcup_{i=1}^N I_i = [0, T]$ and $x$ is monotonic on each $I_i$. $V^1([0,+\ns);\R)$ denotes the subset of functions from $D([0,+\ns);\R)$ with finite total variation on any compact subset of $[0,+\ns)$.

\subsection{Local times}
Let $(E, \cal E)$ be a measurable space. We start with a definition of local time for a deterministic Borel path $x:[0, +\ns) \ra E$ and a given Borel measure $\xi$ on $[0, +\ns)$. For each $t \ge 0$ we  define the corresponding occupation measure $\mu_t$ on $(E, \cal E)$ as 
\begin{equation} \label{eq_local_det}
\forall \;\Gamma\in {\cal E}, \qquad \mu_t(\Gamma) = \int_{[0,t]} {\mathbf 1}_{\Gamma}(x_s) \xi(\dd s) = \xi\rbr{ \cbr{s\in [0, t] : x_s \in \Gamma} }.
\end{equation}
\begin{dfn} \label{def_local_det}
Let $\nu$ be a measure on $(E, \cal E)$ and $\xi$ be a Borel measure on $[0, +\ns)$. We say that $x$ possesses a \emph{local time} $L$ relative to the measures $\nu$ and $\xi$ if for each $t$, $\mu_t \ll \nu$, where $\mu_t$ is the occupation measure defined by \eqref{eq_local_det}. This local time is a function of the space variable $y \in E$ and the time variable $t \in [0, +\ns)$, and is defined as the Radon-Nikodym derivative $L_t^y :=(\dd \mu_t / \dd \nu) (y)$.
\end{dfn}
Immediately from \eqref{eq_local_det} and Definition \ref{def_local_det} we get that for any non-negative Borel function $f: E \ra \R$
\begin{equation} \label{occupation_times_det}
\int_{[0,t]} f\rbr{x_s} \xi(\dd s) = \int_{E} f(y) \mu_t (\dd y) = \int_{E} f(y)  L_t^y \nu (\dd y). 
\end{equation}
The relation \eqref{occupation_times_det} is referred to as \emph{occupation times formula}.

For our purposes, we will need a sufficiently general definition of local time of a stochastic process.
Let $(\Omega, {\cal F}, \P)$ be a probability space and let $X_t$, $t \in [0, +\ns)$, be a stochastic process with measurable state space $(E, \cal E)$, such that the mapping $[0, +\ns)  \times \Omega \ni (t, \omega) \mapsto X_t(\omega)$ is measurable relative to ${\cal B}([0, +\ns)) \otimes {\cal F}$ and $\cal E$ (by ${\cal B}([0, +\ns))$ we denote the Borel $\sigma$-field of $[0, +\ns)$). Let $\Xi$ be a mapping from $\Omega$ to the set of Borel measures on $[0, +\ns)$, which means that for each $\omega \in \Omega$, $\xi = \Xi(\omega)$ is some Borel measure on $[0, +\ns)$. Each trajectory $s \mapsto X_s(\omega)$ is a Borel function and for each $(t, \omega) \in [0, +\ns)  \times \Omega$ we may define a corresponding occupation measure $\mu_t$ on $(E, \cal E)$ as in \eqref{eq_local_det} with $x = X(\omega)$ and $\xi = \Xi(\omega)$.
\begin{dfn} \label{localt_process}
Let $\nu$ be a measure on $(E, \cal E)$. We say that $X$ possesses a \emph{local time} relative to the measure $\nu$ and the mapping $\Xi$ if for each $t \in [0, +\ns)$, $\mu_t \ll \nu$ a.s. (i.e., $\P$-a.s.), where $\mu_t$ is the occupation measure defined as in \eqref{eq_local_det} with $x = X(\omega)$ and $\xi = \Xi(\omega)$, $\omega \in \Omega$.
\end{dfn}
\begin{remark}
Usually, for example when $X$ is a Markov process, the measure $\xi = \Xi(\omega)$ appearing in our definition is equal to the Lebesgue measure (natural clock), and does not depend on $\omega$, see for example \cite{GemanHorowitz80}. However, if one considers the \emph{semimartingale local time} then the measure $\xi = \Xi(\omega)$ is the Lebesgue-Stieltjes measure associated with the continuous part of the quadratic variation of $X$ (business clock), thus may depend on $\omega$.
\end{remark}

\subsection{Truncated variation and numbers of interval crossings} \label{tv_and_crossings}

Now, we define \emph{numbers of interval (up- and down-) crossings} by $x \in D([0,+\ns);\R)$ on the time interval $[s, t]$. 
Let $y\in \R, c > 0$, define $\sigma_{0}^{c}=s$ and for $n \in \N_0$ 
\[
\rho_{n}^{c}:=\inf\left\{ u \in \sbr{ \sigma_{n}^{c}, t}: x_u\ge y+c/2\right\} ,
\]
\[
\sigma_{n+1}^{c}:=\inf\left\{ u \in \sbr{ \rho_{n}^{c}, t} : x_u < y-c/2\right\},
\]
where and throughout the paper we apply conventions: $\inf \emptyset := +\ns$ and $\sbr{+\ns, t} := \emptyset$ for any $t \in [0, +\ns]$.
\begin{dfn}\label{defd} 
The \emph{number of downcrossings, upcrossings and crossings} by $x$ the interval $[y-c/2, y+c/2]$ on the time interval $[s,t]$ are defined as respectively
\begin{align*}
\Dcrossemph{y,c}{x}{[s,t]} :=\max&\left\{ n:\sigma_{n}^{c}\leq t\right\} , \qquad
\Ucrossemph{y,c}{x}{[s,t]} := \Dcrossemph{-y,c}{-x}{[s,t]} \qquad \text{and}\\
&\crossemph{y,c}{x}{[s,t]}:= \Ucrossemph{y,c}{x}{[s,t]} + \Dcrossemph{y,c}{x}{[s,t]}.
\end{align*}

\end{dfn}

The numbers of interval crossings by a real c\`adl\`ag path $x:[0, +\ns) \ra \R$ are closely related to {truncated variation}. 
The relation between the truncated variation and interval crossings is the following. If $x:[0, +\ns) \ra \R$  is a regulated path, then the following formula holds:
\begin{equation} \label{BInd_gen}
\TTV x{\left[s,t\right]}{c} = \int_{\R} \cross{y,c}{x}{[s,t]} \dd y.
\end{equation}
Equation \eqref{BInd_gen} is a generalization of the Banach Indicatrix Theorem, see \cite[Theorem 1]{LochowskiColloquium:2017}.
From \eqref{tvc_def-eq} and the fact that for any $0 \le s < t < u$, $\pi \in \Pi(s,t)$, $\rho \in \Pi(t,u)$, $\pi \cup \rho \in \Pi(s,u)$, it also follows that the truncated variation is a superadditive functional of intervals, that is:
\begin{equation} \label{superadd}
\TTV x{\left[s,u\right]}{c} \ge \TTV x{\left[s,t\right]}{c} + \TTV x{\left[t,u\right]}{c}.
\end{equation}
On the other hand, since for any $0 \le s < t' \le t \le t'' < u$ one has 
\[
\max \rbr{\left| x_{t''}-x_{t'} \right| -c,0} \le \max \rbr{\left| x_{t}-x_{t'} \right| -c,0} + \max \rbr{\left| x_{t''}-x_{t} \right| -c,0} +c,
\]
we have 
\begin{equation} \label{superaddd}
\TTV x{\left[s,u\right]}{c} \le \TTV x{\left[s,t\right]}{c} + \TTV x{\left[t,u\right]}{c} + c.
\end{equation}
\subsection{Truncated variation and variations of a continuous function along Lebesgue partitions} \label{psivar}
When $x$ is continuous, relation \eqref{BInd_gen} allows easily to relate the truncated variation with the Lebesgue partitions of the half line $[0, +\ns]$. 
The Lebesgue partition $\pi^{c,r}$ (sometimes, for typographical reasons, we will also denote it by $\pi(c,r)$) of the interval $[0, +\ns]$ corresponding to $x \in C([0, +\ns); \R)$ and the grid $$c \cdot \Z +r = \cbr{ z \in \R: \exists p \in \Z, z = p\cdot c +r},$$ where $c>0$ and $r \in [0, c)$, is the following family of intervals 
$$\pi^{c,r} = \cbr{\sbr{\tau_k^{c,r}, \tau_{k+1}^{c,r}}, k \in \N_{0}},$$ where $\tau_0^{c,r} = 0$, and for $k \in \N$
\[
\tau^{c,r}_{k+1} :=
\inf\cbr{ t > \tau^{c,r}_{k}:  x_t \in  \rbr{c \cdot \Z +r } \setminus \cbr{x_{\tau^{c,r}_{k}}} }.
\]
For $t\ge 0$ let us define
$k^{c,r}(t) := \max\cbr{k \in \N_0 : \tau^{c,r}_{k} \le t}$. Notice that 
\begin{equation} \label{eesstt1}
\sum_{p \in \Z} \cross{p \cdot c + c/2 +r , c }{x}{[0,t]} \le k^{c,r}(t).
\end{equation}
This follows from the fact that each crossing of the interval $[p\cdot c + r, (p+1)\cdot c + r]$, $p \in \Z$, corresponds to at least one new time $\tau^{c,r}_{k}$. On the other hand, for any $c_0 \in (0, c)$,
\begin{equation} \label{eesstt2}
\sum_{p \in \Z} \cross{p \cdot c + c/2 +r , c_0 }{x}{[0,t]} \ge k^{c,r}(t) - 1
\end{equation}
since each new time $\tau^{c,r}_{k}$, $k \ge 2$, corresponds to at least one crossing of the interval $[p\cdot c + r + (c - c_0)/2, (p+1)\cdot c + r - (c - c_0)/2]$ for some $p \in \Z$.
Now, using \eqref{BInd_gen}, \eqref{eesstt1} and \eqref{eesstt2} we estimate 
\begin{align}
\TTV x{\left[0,t\right]}{c}  & = \int_{\R} \cross{y , c }{x}{[0,t]} \dd y =  \sum_{p \in \Z} \int_{[0, c)} \cross{p \cdot c + c/2 +r , c }{x}{[0,t]} \dd r  \le  \int_{[0, c)}  k^{c,r}(t) \dd r  \qquad\qquad \text{and}\label{eeesstt1}
\end{align}

\begin{align*}
\TTV x{\left[0,t\right]}{c_0}  & =  \int_{\R} \cross{y , c_0 }{x}{[0,t]}  \dd y =  \sum_{p \in \Z} \int_{[0, c_0)} \cross{p \cdot c_0 + c_0/2 +r , c_0 }{x}{[0,t]} \dd r 
\\
& = \sum_{p \in \Z} \int_{[0, c)} \cross{p \cdot c + c/2 +r , c_0 }{x}{[0,t]} \dd r \ge  \int_{[0, c)} \rbr{ k^{c,r}(t) - 1 } \dd r.
\end{align*}
Since the function $(0, +\ns) \ni c \mapsto \TTV x{\left[0,t\right]}{c}$ is continuous (see \cite[Lemma 21]{LochowskiGhomrasniMMAS:2015}), we have 
\begin{align}
\TTV x{\left[0,t\right]}{c}  & = \lim_{c_0 \ra c-} \TTV x{\left[0,t\right]}{c_0} \ge  \int_{[0, c)} \rbr{ k^{c,r}(t) - 1 } \dd r. \label{eeesstt2}
\end{align}

Let now $\psi:[0, +\ns) \ra [0, +\ns)$ be a continuous at $0$, non-decreasing function, starting from $0$, that is $\psi(0) = 0$. Last two relations allow to establish a relationship between $\psi$-variations of $x$ along the Lebesgue partitions $\pi^{c,r}$, defined as
\begin{align*}
V^{c,r}_{\psi}(x, [0,t]) := \sum_{[u, v] \in \pi^{c,r}} \psi\rbr{\left| x_{v \wedge t} - x_{u \wedge t} \right|} =  \sum_{k=1}^{+\ns} \psi\rbr{\left| x_{\tau^{c,r}_{k} \wedge t} - x_{\tau^{c,r}_{k-1} \wedge t} \right|}, \quad   c > 0,\; r\in [0,c),
\end{align*}
and the truncated variation.
Since $\left| x_{\tau^{c,r}_{k}} - x_{\tau^{c,r}_{k-1}} \right| = c$ for $k$ satisfying $k \ge 2$ and $k \le k^{c,r}(t)$, and since $\left| x_{\tau^{c,r}_{1} \wedge t} - x_{0} \right| \le c$,  $\left| x_{t} - x_{\tau^{c,r}_{k^{c,r}(t)}} \right| < c$ and $\tau^{c,r}_{k} \wedge t = t$ for $ k > k^{c,r}(t)$ it is easy to see that the following estimates hold
\begin{align}
V^{c,r}_{\psi}(x, [0,t])   
\in \sbr{ \rbr{k^{c,r}(t) - 1} \psi(c) , \rbr{k^{c,r}(t) + 1} \psi(c) } \label{eeesstt3}.
\end{align}
Using \eqref{eeesstt1}-\eqref{eeesstt3} we have  
\begin{align*}  
&  \psi(c)\TTV x{[0, t]}c - \psi(c)c  \le  \int_{[0, c)} \psi(c) \rbr{ k^{c,r}(t) - 1 } \dd r \le   \int_{[0, c)}  V^{c,r}_{\psi}(x, [0,t])  \dd r   \\ &   \le \int_{[0, c)} \psi(c)  \rbr{ k^{c,r}(t) + 1 }  \dd r  \le   \psi(c)\TTV x{[0, t]}c + 2\psi(c)c. 
\end{align*}
Thus the mean of $\psi$-variations of continuous $x$ along Lebesgue partitions $\pi^{c,r}$, $r \in [0, c)$, is comparable for $c$'s close to $0$ with: 
\begin{align*}  
\frac{ \psi(c)}{c}\TTV x{[0, t]}c - \psi(c) & \le  \frac{1}{c} \int_{[0, c)}  V^{c,r}_{\psi}(x, [0,t])  \dd r  \le \frac{ \psi(c)}{c} \TTV x{[0, t]}c + 2\psi(c). 
\end{align*}
Replacing the variable $r$ with  $\gamma :={r}/{c}$, $\gamma \in [0,1)$, we have
 \begin{align} \label{eeesstt} 
\frac{ \psi(c)}{c}\TTV x{[0, t]}c - \psi(c) & \le  \int_{[0, 1)}  V^{c,\gamma \cdot c}_{\psi}(x, [0,t])  \dd \gamma   \le \frac{ \psi(c)}{c} \TTV x{[0, t]}c + 2\psi(c).  
 \end{align}
In particular, when $\psi(c) = \psi_2(c) = c^2$ we get that $c \cdot \TTV x{[0, t]}c$ is comparable for $c$'s close to $0$ with the mean of  quadratic variations $[x]_t^{\pi\rbr{c,\gamma \cdot c}} := V^{c,\gamma \cdot c}_{\psi_2}(x, [0,t])$ of $x$ along all shifted Lebesgue partitions $\pi\rbr{c,\gamma \cdot c}$, $\gamma \in [0, 1)$:
 \begin{equation*}  
c \cdot \TTV x{[0, t]}c - c \le  \int_{[0, 1)} [x]_t^{\pi\rbr{c, \gamma \cdot c}}  \dd \gamma 
 \le {c} \cdot \TTV x{[0, t]}c + 2c. 
 \end{equation*}
The quantities of the form $\varphi(c)\TTV x{[0, t]}c$ will be essential in the subsequent considerations.

\subsection{Auxiliary results} \label{aux_res}
Recall that $V^1([0,+\ns);\R)$ denotes the subset of \cad functions with finite total variation on any compact subset of $[0,+\ns)$.

For any $x \in V^1([0,+\ns);\R)$ we denote with $\UTV{x}{\dd u}{}$  the Lebesgue-Stieltjes measure  associated to the increasing, right-continuous map $t\mapsto \UTV x{[0, t]}{}$, thus for any $[s,t] \subset [0,+\ns)$
\[
\int_{(s,t]} \UTV{x}{\dd u}{} = \UTV{x}{(s,t]}{}. 
\]
Analogously, one can define $\DTV{x}{\dd u}{}$ and $\TTV{x}{\dd u}{}$.

Using numbers of interval (up-, down-) crossings we define {\em level (up-, down-) crossings}.
  \begin{dfn}\label{defd1} Let $y\in\R$. The number of times that the function $x$ {\em upcrosses} the level $y$ during the time interval $[s, t]$ is defined as 
	 \[
	 	\Ucrossemph{y}{x}{[s,t]} := \lim_{c\to 0+}\Ucrossemph{y,c}{x}{[s,t]}\in \N_0 \cup\cbr{+\ns}.
	 \]
	 Analogously we define the number of {\em downcrosses} as
	 \[
	 	\Dcrossemph{y}{x}{[s,t]} := \lim_{c\to 0+}\Dcrossemph{y,c}{x}{[s,t]}\in \N_0 \cup\cbr{+\ns}.
	 \]
\end{dfn}
\noindent Note that $\Ucross{y,c}{x}{[s,t]}$ and $\Dcross{y,c}{x}{[s,t]} \in \mathbb{N}_0$ are increasing functions of $c$, so the above limits always exist.


Recall that $V^0([0,+\ns);\R)$ denotes the piecewise monotonic \cad functions. The following lemma is required to prove Theorem \ref{thm:meta1}.
\begin{lem} \label{Ban_Viatli_extended} 
Let $x\in V^0([0,+\infty); \R)$, $t>0$ and $g:\R\mapsto \R$ be locally bounded and Borel--measurable. Then
\begin{align}
\int_{\R}g(z)\Ucrossemph{z}{x}{[0,t]}\dd z & =\int_{(0,t]}g(x_{s-})\emph{UTV}(x,\dd  s) +\sum_{0<s\leq t,\Delta x_s>0}\int_{x_{s-}}^{x_s}[g(z)-g(x_{s-})]\dd z,\label{eq:Ban_Vitali_utv}
\end{align}
and analogously 
\begin{align}
\int_{\R}g(z)\Dcrossemph{z}{x}{[0,t]}\dd z & =\int_{(0,t]}g(x_{s-})\emph{DTV}(x,\dd  s)+\sum_{0<s\leq t,\Delta x_s<0}\int_{x_{s}}^{x_{s-}}[g(z)-g(x_{s-})]\dd z.\label{eq:Ban_Vitali_dtv}
\end{align}
\end{lem}
\begin{remark} It is possible to generalize the above lemma for $x \in  V^1([0,+\infty); \R)$, see \cite[Theorem 2.8]{Hove:2024}. For reader's convenience we present its proof for $x \in  V^0([0,+\infty); \R)$.
Lemma \ref{Ban_Viatli_extended}  may also be derived from the change of variable formulas for \cad functions with finite variation proven in \cite[Proposition 1 and Proposition 2]{BertoinYor:2014}. The main difference between \cite{BertoinYor:2014} and Lemma \ref{Ban_Viatli_extended} is such that the autors of \cite{BertoinYor:2014} consider the continuous part of $\TTVemph{x}{\dd u}{}$, hence their formulas have the term $g\rbr{x_s}$, whereas Lemma \ref{Ban_Viatli_extended} uses the term $g\rbr{x_{s-}}$.
\end{remark}

\noindent{\bf Proof of Lemma \ref{Ban_Viatli_extended}:} 
By the assumption $x \in V^0([0,+\8);\R)$ there exists a \emph{finite} sequence of intervals $\{I_i\}_{i=1}^N$, $N \in \N $, which are mutually disjoint, $\bigcup_{i=1}^N I_i = [0,t]$ and the function $x$ is monotone on any of $I_i$. Since $x$ is \cad we may and will assume that these intervals, except the last one containing $t$, are of the form $I_i = \left[t_i, t_{i+1}\right)$, where $t_i < t_{i+1}$. Moreover, we will assume that they are the largest intervals possible on which $x$ is monotone. Thus, if for some $i>1$, $x$ is non-decreasing on $ \left[t_{i-1}, t_{i}\right)$ and non-decreasing on $ \left[t_i, t_{i+1}\right)$ then there need to be a downward jump at the time $t_i$. Similarly, if for some $i>1$, $x$ is non-increasing on $ \left[t_{i-1}, t_{i}\right)$ and non-increasing on $ \left[t_i, t_{i+1}\right)$  then there need to be an upward jump at the time $t_i$. Let $S$ denote the set of such times $s$. 

Let us define $R = \bigcup_{i=1}^N \cbr{x_{t_i}, x_{t_{i+1}}}$. Given a level $z \in \R \setminus R$, with any interval $I_i$ or any $s \in S$ we can have at most one associated upcrossing. This happens when $x_{t_i} < x_{t_{i+1}}$ or $x_{s-} < x_s$ and
	\[
		z\in \rbr{x_{t_i}, x_{t_{i+1}}} \text{ or } z\in \rbr{x_{s-}, x_s}.
	\]  
	Clearly, neither of the intervals where the function is decreasing induces any upcrossing. Let $J$ be the subset of the indices where the function is increasing and $T$ be the subset of $s \in S$ such that $x$ has an upward jump at $s$.  Then we have
	\[
		\Ucross{z}{x}{[0,t]} = \sum_{i\in J} \mathbf{1}_{ \rbr{x_{t_i}, x_{t_{i+1}}}  } (z) + \sum_{s\in T} \mathbf{1}_{ \rbr{x_{s-}, x_s} } (z) .
	\]
	Consequently
	\[
		\int_{\R} g(z)\Ucross{z}{x}{t} \dd{z} = \int_{\R \setminus R} g(z)\Ucross{z}{x}{t} \dd{z} = \sum_{i\in J} \int_{x_{t_i}}^{x_{t_{i+1}}} g(z) \dd{z} + \sum_{s\in T} \int_{x_{s-}}^{x_{s}} g(z) \dd{z}.
	\]
	We are now to deal with these integrals. To this end, we apply the classical idea of opening temporal windows at times of jumps. We may impose that the sum of the lengths of these windows is finite and consider interpolated continuous $\tilde x$. We have
\[
\int_{x_{t_i}}^{x_{t_{i+1}}} g(z) \dd{z} = \int_{\tilde{x}_u}^{\tilde{x}_{u'}}g(z)\dd z \text{ and }
\int_{x_{s-}}^{x_{s}} g(z) \dd{z} = \int_{\tilde{x}_v}^{\tilde{x}_{v'}}g(z)\dd z
\]
for properly defined $u,u', v, v'$. Then we apply the change of variable for the Riemann-Stieltjes integral 
\[
\int_{\tilde{x}_u}^{\tilde{x}_{u'}}g(z)\dd z  = \int_{u}^{u'} g{\rbr{\tilde{x}_s}} \dd \tilde{x}_{s} \text{ and } \int_{\tilde{x}_v}^{\tilde{x}_{v'}}g(z)\dd z  = \int_{v}^{v'} g{\rbr{\tilde{x}_s}} \dd \tilde{x}_{s}
\]
Clearly, for properly defined $\tilde{t}$ we have
\begin{equation} \label{change_of_var}
\int_{\R} g(z)\Ucross{z}{x}{t} \dd{z} = \int_{0}^{\tilde{t}}  g{\rbr{\tilde{x}_s}} \dd \tilde{x}_{s}.
\end{equation}
Let us consider the decomposition of the measure $\text{UTV}(x,\dd s)$ into the continuous part $\mu^c$ and the atomic part $\mu^a.$ Moreover, let $W$ be the set of the temporal windows. We obtain
\[
	 \int_{0}^{\tilde{t}} \mathbf{1}_{s\notin W} g{\rbr{\tilde{x}_s}} \dd \tilde{x}_{s} = \int_0^t g(x_{s-}) \dd \mu^c\rbr{s},
\]
moreover
\[
	\int_{0}^{\tilde{t}} \mathbf{1}_{s\in W} g{\rbr{\tilde{x}_s}} \dd \tilde{x}_{s}  = \sum_{0 < s \le t, \Delta x_s>0} \int_{x_{s-}}^{x_s} g(z) \dd{z}.
\]
Using this and (\ref{change_of_var}) we obtain \eqref{eq:Ban_Vitali_utv} by simple calculations. Similarly one can prove (\ref{eq:Ban_Vitali_dtv}).
\hfill $\square$

To state the next result, we will need some properties of the \emph{Skorohod problem on $[-c/2, c/2]$, $c >0$,} proven in \cite{LochowskiGhomrasniMMAS:2015}, also see \cite{burdzy:2009}. Let us first recall the definition of the  Skorohod map. 
\begin{dfn}
Let $x, \alpha,\beta \in D([0,\8);\R)$ and $\phi_0\in\mathbb{R}.$ A pair
of functions $\left(\phi^{x},\eta^{x}\right)\in D[0;+\infty)\times V^1([0,+\ns);\R)$
is said to be a solution of the \emph{Skorohod problem on $\left[\alpha, \beta\right]$
with the starting condition $\phi^{x}(0)=\phi_0$ for $x$} if the following
conditions are satisfied: 
\begin{enumerate}
\item 
for every $t\geq0,$ $\phi^{x}\left(t\right)=x\left(t\right)+\eta^{x}\left(t\right)\in\left[\alpha\left(t\right), \beta\left(t\right)\right];$ 
\item 
$\eta^{x}=\eta_{d}^{x}-\eta_{u}^{x},$ where $\eta_{d}^{x},\eta_{u}^{x}$ are non-decreasing \cad functions 
and the corresponding measures $\mathrm{d}\eta_{d}^{x},$ $\mathrm{d}\eta_{u}^{x}$
are carried by $\left\{ t\geq 0:\phi^{x}(t)=\alpha(t)\right\} $ and $\left\{ t\geq 0:\phi^{x}(t)=\beta(t)\right\} $
respectively;
\item 
$\phi^{x}(0)=\phi_0.$
\end{enumerate}
\end{dfn}
\noindent The Skorohod problem on $\left[\alpha, \beta\right]$ has always a solution whenever 
\[
\inf_{t \ge 0} \cbr{ \beta(t) - \alpha(t) } >0 \text{ and } \phi_0 \in [\alpha(0), \beta(0)],
\]
see \cite[Proposition 2.7]{LochowskiGhomrasniMMAS:2015}.
In such a case $\eta^x \in V^0([0,+\ns);\R)$ and the function $-\eta^x$ is called the \emph{regularisation of $x$ obtained via the Skorohod map on $[\alpha, \beta]$} (with the starting condition $\phi_0 =x(0) + \eta^x(0) $).
In the case when $\alpha \equiv -c/2$, $\beta \equiv c/2$, $c>0$, the problem is called the Skorohod problem on $[-c/2, c/2]$.

Using Lemma \ref{Ban_Viatli_extended} we will prove Theorem \ref{thm:meta1}, relating the integrated numbers of interval crossings of a function $x\in D([0,\8);\R)$ with the asymptotics of the upward and downward truncated variations $\UTV x{t}{c},$  $\DTV x{t}{c},$ $t \ge 0,$ as $c \ra 0+$.
\\{\bf Proof of Theorem \ref{thm:meta1}:} 
Conditions  \eqref{conv_tv0}, \eqref{conv_utv} and \eqref{conv_dtv} are equivalent. This follows from equalities (see \cite[Sect. 2]{LochowskiGhomrasniMMAS:2015}): 
\[
\TTV x{[0,t]}c = \UTV x{[0,t]}c+ \DTV x{[0,t]}c 
\quad \text{and} \quad  x_t^c - x_0^c = \UTV x{[0,t]}c - \DTV x{[0,t]}c \]
($x^c$ is a regularisation of $x$ obtained via the Skorohod map on $[-c/2, c/2]$ for some properly chosen starting condition).
To get equivalence of \eqref{conv_tv0}, \eqref{conv_utv} and \eqref{conv_dtv}  it is sufficient to notice that the differences 
\[
\varphi(c)\UTV x{[0,t]}c - \varphi(c)\DTV x{[0,t]}c = 2\varphi(c)\UTV x{[0,t]}c - \varphi(c)\TTV x{[0,t]}c  = \varphi(c)\rbr{x_t^c - x_0^c}
\]  
tend to $0$ as $c \ra 0+$.

Assume now that (\ref{conv_utv}) holds. The fact that $\zeta$ is non-negative and non-decreasing follows immediately from the same properties of the function $t \mapsto \UTV x{[0,t]}c,$ $t \ge 0.$ 

Now let us fix $T>0$ and put
\begin{equation*}
 M= \max\left\{\sup_{0\leq t\leq T}|x_t|, \sup_{0\leq t\leq T}|\zeta_t| = \zeta_T \right\}.
\end{equation*}
Using the assumption that $g$ is continuous, for fixed $\varepsilon>0$ we can find a number $\delta>0$ such that
\begin{equation}\label{delta_estimate}
\sup_{z,z'\in[-M-1;M+1],|z-z'|\leq\delta}|g(z)-g(z')| \leq \frac{\varepsilon}{M}.
\end{equation}
Let us also define
\begin{equation}
N = \#\left\{ t\in(0;T]:|x_t-x_{t-}|>\delta\right\}.
\end{equation}
Let $t\in[0,T].$ We will use the regularisation $x^c$ of $x$ obtained via the Skorokhod map on $[-c/2, c/2]$. Since for any $z \in \R$ and any $t >0$, $\Ucross{z,c}{x}{[0,t]} = \Ucross{z}{x^c}{[0,t]} \pm 1$ (\cite[Lemma~3.3 and~3.4]{LochowskiGhomrasni:2014}) and $\Ucross{z,c}{x}{[0,t]} = \Ucross{z}{x^c}{[0,t]}=0$ when $z<\inf_{s\in[0,t]}x_s-c/2$ or $z>\sup_{s\in[0,t]}x_s+c/2$ to prove convergence \eqref{weak_local_tv} when $\mathrm{n}^{z, c}(x,[0, t])$ is replaced by $2\mathrm{u}^{z, c}(x,[0, t])$, it is sufficient to prove the convergence  
\[
\varphi(c)\int_{\R}\Ucross{z}{x^c}{[0,\cdot]}g(z)\dd z \rightarrow \frac{1}{2} \int_{0}^{\cdot}g(x_{t-})\dd \zeta_t \text{ as } c \ra 0+.
\] 
To prove this convergence we will use Lemma \ref{Ban_Viatli_extended}. From the properties of the Skorokhod map we have that $x^c$ is piecewise monotonic (see \cite[Proposition 2.7 and formula (2.4)]{LochowskiGhomrasni:2014}). Using Lemma \ref{Ban_Viatli_extended} we write
\begin{align}
\int_{\R}g(z) \Ucross{z}{x^{c}}{[0,t]} \dd{z} = & \int_{0}^{t}g(x^c_{s-})\UTV{x^c}{\dd  s}{} +\sum_{0<s\leq t,\Delta x^c_s>\delta}\int_{x^c_{s-}}^{x^c_s}[g(z)-g(x^c_{s-})]\dd z\nonumber \\
 & +\sum_{0<s\leq t, 0 <\Delta x^c_s \le \delta}\int_{x^c_{s-}}^{x^c_s}[g(z)-g(x^c_{s-})]\dd z.\label{eq:trzyyyy}
\end{align}
Let us denote the consecutive summands on the right side of equation (\ref{eq:trzyyyy}) by $R_1(c)$, $R_2(c)$ and $R_3(c)$ respectively. 

First, for $c>0$ we estimate
\begin{equation}
    \varphi(c)R_2(c)  \leq\varphi(c)\cdot2N(2M+c)\rbr{2\sup_{z\in[-M-c/2;M+c/2]}|g(z)|}\ra 0 \quad \text{as }c\ra 0+.\label{eq:estimate0}
\end{equation}
Using (\ref{delta_estimate}) and the fact that $\UTV{x^c}{[0, t]}{} \le \UTV{x}{[0, t]}{c} + c$ (see \cite[Proposition 2.9]{LochowskiGhomrasni:2014}) we estimate 
\begin{align}
\varphi(c)R_3(c)  & \leq\varphi(c)\sum_{0<s\leq t,0<\Delta x^{c}_s\leq\delta}|\Delta x^{c}_s|\sup_{z,z'\in[-M-c/2;M+c/2],|z-z'|\leq\delta}|g(z)-g(z')|\nonumber \\
 & \leq\varphi(c)\UTV{x^{c}}{[0,t]}{}\frac{\varepsilon}{M}
 \leq\varphi(c)(\UTV x{[0,t]}c +c) \frac{\varepsilon}{M} \le 2 \varepsilon \label{eq:estimate1}
\end{align}
for  $c$ small enough such that $\varphi(c)\UTV x{[0,T]}c + \varphi(c) c \le 2 M.$

Now we are left with $R_1$. We fix $K = 1,2,\ldots $ and define 
\[
I_{i}=\left[-M+2M\frac{i-1}{K};-M+2M\frac{i}{K}\right],\quad i=\{1,2,...,K\}.
\]
Further we define the following sequence $\cbr{\upsilon_k}_{k\geq 0}$ of times. Let $\upsilon_{0}=0$ and 
\begin{equation*}
\upsilon_{k}=\inf\left\{ t>\upsilon_{k-1}:x_t\in I_{i}\mbox{ for some }i=1,2,...,K\mbox{ such that }x_{\upsilon_{k-1}}\notin I_{i}\right\} .
\end{equation*}
We have 
\begin{align}
R_1(c)= & \sum_{k=0}^{+\infty}\int_{(\upsilon_{k}\wedge t,\upsilon_{k+1}\wedge t]}g(x^{c}_{s-}) \UTV{x^{c}}{\dd s}
\nonumber \\
 &=  \sum_{k=0}^{+\infty}\int_{(\upsilon_{k}\wedge t,\upsilon_{k+1}\wedge t]}\left\{ g(x^{c}_{s-})-g(x_{\upsilon_{k}\wedge t})\right\} \UTV{x^{c}}{\dd s} \; + \sum_{k=0}^{+\infty}g(x_{\upsilon_{k}\wedge t})\UTV{x^{c}}{(\upsilon_{k}\wedge t,\upsilon_{k+1}\wedge t]}{},\label{eq:estimate2}
\end{align}
where we denote $\UTV{x^{c}}{(\upsilon_{k}\wedge t,\upsilon_{k+1}\wedge t]}{} = \UTV{x^{c}}{[0,\upsilon_{k+1}\wedge t]}{} - \UTV{x^{c}}{[0,\upsilon_{k}\wedge t]}{} $. 
Using (\ref{delta_estimate}) it is easy to estimate for $c \le 2\delta$ and $2M/K \le \delta$ the first summand in (\ref{eq:estimate2}) multiplied by $\varphi(c):$
\begin{align}
 & \varphi\left(c\right)\sum_{k=0}^{+\infty}\int_{(\upsilon_{k}\wedge t,\upsilon_{k+1}\wedge t]}\left|g(x^{c}_{s-})-g(x_{s-})+g(x_{s-})-g(x_{\upsilon_{k}\wedge t})\right| \UTV{x^{c}}{\dd s}\nonumber \\
 & \leq\varphi\left(c\right)\sum_{k=0}^{+\infty}\frac{2 \varepsilon}{M} \UTV{x^{c}}{\left(\upsilon_{k}\wedge t,\upsilon_{k+1}\wedge t\right]}{} \le  \varphi\left(c\right) \UTV{x^{c}}{[0,t]}{} \frac{2 \varepsilon}{M}  \leq \varphi\left(c\right) \rbr{\UTV{x}{[0,t]}{c} +c } \frac{2 \varepsilon}{M} \le 4 \varepsilon \label{eq:estim2}
\end{align}
for  $c$ small enough such that $\varphi(c)\UTV x{[0,T]}c + \varphi(c) c \le 2 M.$
Now we investigate the second summand of $R_1$ multiplied by $\varphi\left(c\right)$. By \eqref{conv_utv}:
\begin{align}
\sum_{k=0}^{+\infty}g(x_{\upsilon_{k}\wedge t})\varphi(c) \text{UTV}(x^{c},(\upsilon_{k}\wedge t,\upsilon_{k+1}\wedge t])
\; \rightarrow \; \frac{1}{2}\sum_{k=0}^{+\infty}g(x_{\upsilon_{k}\wedge t}) \left\{ \zeta_{\upsilon_{k+1}\wedge t}-\zeta_{\upsilon_{k}\wedge t}\right\} \quad \text{ as }c\ra 0+. \label{eq:estim3}
\end{align}
Moreover for $2M/K\le \delta$ 
\begin{align}
& \left|\int_{(0,t]}g(x_{s-})\dd \zeta_{ s}-\sum_{k=0}^{+\infty}g(x_{\upsilon_{k}\wedge t})\right. \left. \left\{ \zeta_{\upsilon_{k+1}\wedge t}-\zeta_{\upsilon_{k}\wedge t} \right\} \right|\nonumber 
 =\left|\sum_{k=0}^{+\infty}\int_{(\upsilon_{k}\wedge t,\upsilon_{k+1}\wedge t]}\left\{ g(x_{s-})-g(x_{\upsilon_{k}\wedge t})\right\} \dd \zeta_s \right|
 \\
 &\quad \quad \quad \leq\sum_{k=0}^{+\infty}\int_{(\upsilon_{k}\wedge t,\upsilon_{k+1}\wedge t]}\frac{\varepsilon}{M} \dd \zeta_s= \frac{\varepsilon}{M} \zeta_t \le \frac{\varepsilon}{M} \zeta_T \le \varepsilon. \label{eq:estim4}
\end{align}
From (\ref{eq:trzyyyy})-(\ref{eq:estim4})
we get 
\[
\sup_{0 \le t \le T} \left|\varphi(c) \int_{\R}g(z) \Ucross {z,c}x{[0,t]} \dd z-\frac{1}{2} \int_{(0,t]}g(x_{s-}) \dd \zeta_s \right|\leq 10 \varepsilon
\]
for any $c>0$ sufficiently close to $0$. 

When $\zeta$ is continuous, since the functions $t \mapsto \varphi(c) \UTV{x^{c}}{[0,t]}{}$ are non-decreasing, the convergence $\varphi(c) \UTV{x^{c}}{[0,t]}{} \ra \zeta_t$ as $c \ra 0+$ is uniform on $[0,T]$. From this and the fact that there is only a finite number of $k \in \N_0$ for which 
$\upsilon_k \le T$, we get that the convergence in (\ref{eq:estim3}) is uniform on $[0,T]$ and, since all other estimates were uniform in $t$, we get that the convergence in \eqref{weak_local_tv},  when $\mathrm{n}^{z, c}(x,[0, t])$ is replaced by $2\mathrm{u}^{z, c}(x,[0, t])$, is uniform on $[0,T]$ as well. 

The proof of the convergence of $\int_{\R}g(z)\Dcross {z,c}x{[0,t]}\dd z$ under assumption \eqref{conv_dtv} is analogous. The convergence of $\int_{\R}g(z)\cross {z,c}x{[0,t]}\dd z$ follows from the relation $\cross {z,c}x{[0,t]} = \Ucross {z,c}x{[0,t]} + \Dcross {z,c}x{[0,t]}$ and the convergences of $\int_{\R}g(z)\Ucross {z,c}x{[0,t]}\dd z$  and $\int_{\R}g(z)\Dcross {z,c}x{[0,t]}\dd z$.

\hfill $\square$
\begin{remark} \label{subsequence}
Theorem \ref{thm:meta1} remains valid when the convergence $c \ra 0+$ is replaced by a convergence along some sequence $c_n \ra 0+$ as $n \ra +\ns.$ To see this it is enough to replace everywhere in the proof $c$ by $c_n$ and $c \ra 0+$ by $n \ra +\ns$.
\end{remark}

\subsection{Examples of application of Theorem \ref{thm:meta1}} \label{subsect:examples}

\textbf{Example 1.} 
Let $\rbr{t_n}_{n\in \N}$ be a strictly increasing sequence ($t_n < t_{n+1}$, $n \in \N_0$) such that $t_0 = 0$ and $\lim_{n \ra +\ns} t_n = 1$. Next, let $\rbr{a_n}_{n\in \N_0}$ be a non-increasing sequence of positive reals such that $\lim_{n \ra +\ns} a_n = 0$ and 
\begin{equation} \label{calc0}
\sum_{n \in \N_0} a_n = +\ns.
 \end{equation}
Define a continuous ''zigzag'' function $x^1:[0,1] \ra \R$,
\[
x^1(t) = 
\begin{cases}
\frac{a_n}{t_{2n+1}-t_{2n}}\rbr{t-t_{2n}} & \text{ if } t \in \left[t_{2n}, t_{2n+1} \right), n\in \N_0;
\\
\frac{a_n}{t_{2n+2}-t_{2n+1}}\rbr{t_{2n+2}-t} & \text{ if } t \in \left[t_{2n+1}, t_{2n+2} \right), n\in \N_0;
\\
0 & \text{ if } t=1.
\end{cases}
\]
On each time interval $\sbr{t_{2n}, t_{2n+1}}$, $n \in \N_0$, there is exactly one upcrossing the value interval $[y-c/2, y+c/2]$ iff $c/2 \le y < a_n - c/2$ and on each time interval $\sbr{t_{2n+1}, t_{2n+2}}$, $n \in \N_0$, there is exactly one downcrossing the value interval $[y-c/2, y+c/2]$ iff $c/2 < y \le a_n - c/2$.  Using this observation and \eqref{BInd_gen} we calculate 
\begin{equation} \label{calc1}
\TTV {x^1}{[0,1]}c = 2 \sum_{n\in \N_0} \max\rbr{a_n-c,0}
\end{equation}
and for $t \in [0,1)$ we estimate
\begin{equation} \label{calc2}
 \TTV {x^1}{[0,t]}c \le 2a_0 \rbr{ \max \cbr{n\in \N_0: t_{n} \le t} + 1}. 
\end{equation}

Let us define $\varphi:(0,+\ns) \ra (0, 1]$ as
\[
\varphi(c) := \frac{1}{1 + \TTV {x^1}{[0,1]}c}.
\]
From \eqref{calc0}-\eqref{calc2} we get
\[
\lim_{c \ra 0+} \varphi(c) \TTV {x^1}{[0,t]}c = \zeta^1_t = 
\begin{cases}
0 & \text{ if } t \in [0,1);
\\
1 & \text{ if } t=1
\end{cases}
\]
and by Theorem \ref{thm:meta1} for any continuous $g:\R \ra \R$ we have 
\[
\varphi(c)\int_{\R} \cross{y,c}{x^1}{[0, t]} g(y) \dd y\rightarrow \int_{(0, t]}g(x^1_{s-})\dd \zeta^1_s = 
\begin{cases}
0 & \text{ if } t \in [0,1);
\\
g\rbr{x^1_{1-} } = g(0)& \text{ if } t=1.
\end{cases}
\]

Let us notice that this example also shows that in the case when the limit $\zeta$ appearing in statements of Theorems \ref{thm:meta1} are not continuous, the convergence in \eqref{weak_local_tv} may not be locally uniform. For example, taking $g \equiv 1$, we have that 

\begin{align*}
\varphi(c)\int_{\R} \cross{y,c}{x^1}{[0, t]} \dd y = \varphi(c) \TTV {x^1}{[0,t]}c   \ra 
\begin{cases}
0 & \text{ if } t \in [0,1);
\\
1& \text{ if } t=1.
\end{cases}
\end{align*}
The functions $[0, 1]\ni t \mapsto \varphi(c) \TTV {x^1}{[0,t]}c$ are continuous for any $c >0$ (since $x^1$ is continuous), thus the convergence can not be uniform.

\textbf{Example 2.} 
In this example, we will use the Cantor set. The Cantor is obtained first by removing the middle third of the segment $[0,1]$, that is the interval $I_{0,0} = \rbr{\frac{1}{3}, \frac{2}{3}}$. From two remaining segments again the middle thirds, that is the intervals $I_{1,0} = \rbr{\frac{1}{3^2}, \frac{2}{3^2}}$ and $I_{1,1} = \rbr{\frac{7}{3^2}, \frac{8}{3^2}}$, are removed. Continuing this process, that is removing the middle thirds of the remaining segments, we obtain the Cantor set
$${\mathcal C} = [0,1] \setminus \bigcup_{n \in \N_0} \bigcup_{k=0}^{2^n-1} I_{n,k}.$$
Let $\rbr{b_n}_{n\in \N_0}$ be a non-increasing sequence of positive reals such that $\lim_{n \ra +\ns} b_n = 0$ and
\begin{equation} \label{calc20}
\sum_{n \in \N_0} 2^n b_n = +\ns.
 \end{equation}

Now let us define 
\[
x^2(t) = 
\begin{cases}
2\cdot {3^{n+1}} b_n  \min_{u \in \mathcal C} |t - u| & \text{ if } t \in I_{n,k}, n\in \N_0, k \in \cbr{0, 1, 2, \ldots, 2^n-1};
\\
0 & \text{ if } t \in \mathcal C.
\end{cases}
\]
On each interval $ \bar{I}_{n,k}, n\in \N_0, k \in \cbr{0, 1, 2, \ldots, 2^n-1}$, the function $x^2$ starts from $0$, then increases to $b_n$ and again decreases to $0$.

Similarly as in the previous example, using \eqref{BInd_gen} we calculate 
\begin{equation} \label{calc21}
\TTV {x^2}{[0,1]}c = \sum_{n\in \N_0} \sum_{k=0}^{2^n-1} 2 \max\rbr{b_n-c,0}
 = \sum_{n\in \N_0} 2^{n+1} \max\rbr{b_n-c,0}.
\end{equation}
For $t \in [0,1]$ and $n \in \N_0$ let $k_{n}(t)$ denote the number of segments $I_{n,k}$, $n\in \N_0$, $k \in \cbr{0, 1, 2, \ldots, 2^n-1}$, such that $I_{n,k} \subseteq [0, t]$:
\begin{equation} \label{kndef}
k_{n}(t) = \max \rbr{\cbr{k\in \cbr{1, 2, \ldots, 2^n}:  I_{n,k-1} \subseteq [0, t]} \cup \cbr{0}}.
\end{equation}
We have the following estimate
\begin{align} 
&  \sum_{n\in \N_0} \sum_{k=0}^{k_n(t)-1} 2 \max\rbr{b_n-c,0}  = \sum_{n\in \N_0}  2 k_{n}(t) \max\rbr{b_n-c,0} 
\nonumber \\ & \le \TTV {x^2}{[0,t]}c  \le \sum_{n\in \N_0} \sum_{k=0}^{k_n(t)} 2 \max\rbr{b_n-c,0} =  \sum_{n\in \N_0}  2 ( k_{n}(t) +1) \max\rbr{b_n-c,0}. \label{calc22}
\end{align}
It is not difficult to notice that if $t \in I_{n,k}$ for some $n\in \N_0$ and $k \in \cbr{0, 1, 2, \ldots, 2^n-1}$ then
\begin{equation} \label{calc23}
 \frac{1}{2} k_{n+1}(t) = \frac{1}{2^2} k_{n+2}(t) = \frac{1}{2^3} k_{n+3}(t)  \text{ etc.} \quad \text{and,}
\quad 
1 =  k_{0}(1) = \frac{1}{2} k_{1}(1) = \frac{1}{2^2} k_{2}(1)  \text{ etc.} 
\end{equation} 
This follows easily from the fact that if $I_{n,k} = \rbr{a^{n,k}, b^{n,k}}$ then the ternary expansion of $a^{n,k}$ and $b^{n,k}$ is of the form 
$a^{n,k}=0.d_1d_2 \ldots d_n 1_{(3)} = \sum_{i=1}^n d_i 3^{-i} + {1\cdot{3^{-(n+1)}}}$, $b^{n,k}=0.d_1d_2 \ldots d_n 2_{(3)} = \sum_{i=1}^n d_i 3^{-i} + 2\cdot{3^{-(n+1)}}$, where $d_1,d_2, \ldots,d_n \in \cbr{0,2}$.

Let us define $\varphi:(0,+\ns) \ra (0, 1]$ as
\[ 
\varphi(c) := \frac{1}{1 + \TTV {x^2}{[0,1]}c}.
\]
From \eqref{calc20} and \eqref{calc22}-\eqref{calc23} we get
\begin{equation} \label{functx2}
\lim_{c \ra 0+} \varphi(c) \TTV {x^2}{[0,t]}c = \zeta^2_t = 
\begin{cases}
\frac{k_{n+1}(t)}{2^{n+1}} & \text{ if } t \in I_{n,k};
\\
\lim_{n \ra + \ns} \frac{k_{n+1}(t)}{2^{n+1}} & \text{ if } t \in \mathcal C.
\end{cases}
\end{equation}
Indeed, for $t \in I_{n_0,k}$, $n_0\in \N_0$, $k \in \cbr{0, 1, 2, \ldots, 2^{n_0}-1}$,  denoting $ \zeta^2_t = \frac{1}{2^{n_0+1}} k_{n_0+1}(t)$ and using \eqref{calc22} we get 
\begin{align} 
& 2 \zeta_t^2 \sum_{n=n_0+1}^{+\ns}  2^n \max\rbr{b_n-c,0}  \le \TTV {x^2}{[0,t]}c  \nonumber \\
 & \le  2\sum_{n=0}^{n_0}  2^{n+1} \max\rbr{b_n-c,0} +2 \zeta_t^2 \sum_{n=n_0+1}^{+\ns}  2^n \max\rbr{b_n-c,0}  +  2  \sum_{n=n_0+1}^{+\ns}   \max\rbr{b_n-c,0}.  \label{calc244}
\end{align}
Let us take $N_0 \in \N_0$, $N_0 \ge n_0$, divide inequalities \eqref{calc244} by $2\sum_{n=n_0+1}^{+\ns}  2^n \max\rbr{b_n-c,0} $ and notice that 
\[
\lim_{c \ra 0+} \frac{\sum_{n=0}^{n_0}  2^{n+1} \max\rbr{b_n-c,0}}{\sum_{n=n_0+1}^{+\ns}  2^n \max\rbr{b_n-c,0}} = \sum_{n=0}^{n_0}  \lim_{c \ra 0+}  \frac{2^{n+1} \max\rbr{b_n-c,0}}{\sum_{n=n_0+1}^{+\ns}  2^n \max\rbr{b_n-c,0}}= 0
\]
and
\[
\lim_{c \ra 0+} \sum_{n=n_0+1}^{N_0} \frac{\max\rbr{b_n-c,0}}{\sum_{n=n_0+1}^{+\ns}  2^n \max\rbr{b_n-c,0} } = 0
\]
(this follows from \eqref{calc20}).
Next, we estimate
\begin{align*}  \frac{\sum_{n=n_0+1}^{+\ns}   \max\rbr{b_n-c,0}}{\sum_{n=n_0+1}^{+\ns}  2^n \max\rbr{b_n-c,0}  } & \le \frac{\sum_{n=n_0+1}^{N_0}   \max\rbr{b_n-c,0}}{\sum_{n=n_0+1}^{+\ns}  2^n \max\rbr{b_n-c,0} }+ \frac{\sum_{n=N_0+1}^{+\ns}   \max\rbr{b_n-c,0}}{\sum_{n=N_0+1}^{+\ns}  2^n \max\rbr{b_n-c,0}  } \\
& \le  \sum_{n=n_0+1}^{N_0} \frac{\max\rbr{b_n-c,0}}{\sum_{n=n_0+1}^{+\ns}  2^n \max\rbr{b_n-c,0} }+ \frac{1 }{2^{N_0+1}}.
\end{align*} 
Since $N_0$ may be arbitrarily large, from the last three calculations and \eqref{calc244} we get
\[
\lim_{c \ra 0+} \frac{\TTV {x^2}{[0,t]}c}{2\sum_{n=n_0+1}^{+\ns}  2^n \max\rbr{b_n-c,0}} = \zeta_t^2.
\]
Finally, noticing that 
\[
\lim_{c \ra 0+} \frac{2\sum_{n=n_0+1}^{+\ns}  2^n \max\rbr{b_n-c,0}}{1+ \TTV {x^2}{[0,1]}c} = \lim_{c \ra 0+} \frac{2\sum_{n=n_0+1}^{+\ns}  2^n \max\rbr{b_n-c,0}}{1+ 2\sum_{n=0}^{+\ns}  2^n \max\rbr{b_n-c,0}} =1
\]
we get \eqref{functx2}.

The function $\zeta^2$ is called the Cantor function a.k.a. `devil's staircase' function. For example $\zeta^2_t = 1/2$ for $t \in I_{0,0}$,  $\zeta^2_t = 1/4$ for $t \in I_{1,0}$, $\zeta^2_t = 3/4$ for $t \in I_{1,1}$ etc. $\zeta^2$ is continuous but the measure $\dd \zeta^2_t$ charges only points from the Cantor set $\mathcal C$ and, by Theorem \ref{thm:meta1}, for any continuous $g:\R \ra \R$ we have 
\begin{align*}
\varphi(c)\int_{\R} \cross {y,c}{x^2}{[0,t]}  g(y) \dd y & \rightarrow \int_{(0, t]}g(x^2_{s-})\dd \zeta^2_s \\ 
& \quad = \int_{\mathcal C \cap (0, t]}g(x^2_{s-})\dd \zeta^2_s =  \int_{\mathcal C \cap (0, t]}g(x^2_{s})\dd \zeta^2_s = \int_{\mathcal C \cap (0, t]}g(0)\dd \zeta^2_s = g(0) \zeta^2_t.
\end{align*}

\textbf{Example 3.} As in Example 1, let $\rbr{t_m}_{m\in \N_0}$ be a strictly increasing sequence ($t_m < t_{m+1}$, $n \in \N$) such that $t_0 = 0$ and $\lim_{m \ra +\ns} t_m = 1$ and let $\rbr{a_m}_{m\in \N}$ be a non-increasing sequence of positive reals such that $\lim_{m \ra +\ns} a_n = 0$ and 
\begin{equation} \label{calc30}
\sum_{m\in \N} a_m = +\ns.
 \end{equation}
Let $\rbr{m_n}_{n \in \N_0}$ be a strictly increasing sequence of positive integers such that
\begin{equation} \label{calc300}
\sum_{n=0}^{+\ns} 2^n a_{m_n} < 1
\end{equation}
(we will use this assumption later in Subsect. \ref{subsect46}).
Let 
\[
I_{n,k} =  \rbr{a^{n,k}, b^{n,k}}, \quad n\in \N_0, \quad k \in \cbr{0, 1, 2, \ldots, 2^n-1}
\] 
be the intervals used in the construction of the Cantor set in Example 2.
For $n\in \N_0, k \in \cbr{0, 1, 2, \ldots, 2^n-1}$  we define
\[
t_l^{n,k} =  a^{n,k} + \rbr{b^{n,k} -a^{n,k}} t_l = a^{n,k} + 3^{-n-1} t_l, \quad l\in \N_0,
\]
and a continuous ''zigzag'' function $x^3:[0,1] \ra \R$ such that 
\[
x^3(t) = 
\begin{cases}
\zeta_t^2  + a_{m_n+l} \times \frac{t-t_{2l}^{n,k}}{t_{2l+1}^{n,k}-t_{2l}^{n,k}} & \text{ if } t \in \left[t_{2l}^{n,k}, t_{2l+1}^{n,k} \right)\cap I_{n,k}, \quad l\in \N_0;
\\
\zeta_t^2  + a_{m_n+l} \times \frac{t_{2l+2}^{n,k} - t}{t_{2l+2}^{n,k}-t_{2l+1}^{n,k}} & \text{ if } t \in \left[t_{2l+1}^{n,k}, t_{2l+2}^{n,k} \right)\cap I_{n,k}, \quad l\in \N_0;
\\
\zeta_t^2 & \text{ if } t \in \mathcal C.
\end{cases}
\]

For $t \in [0,1]$ and $n \in \N_0$ let $k_n(t)$ be defined by \eqref{kndef}.
Similarly as in two previous examples, for $t \in I_{n_0,k}$, $c \in \rbr{0, a_{m_{n_0}}}$, $n_0 \in \N_0$, $k \in \cbr{0, 1, 2, \ldots, 2^{n_0}-1}$,  we have the following estimate
\begin{align} 
&  \sum_{n = n_0+1}^{+\ns} \sum_{k=0}^{k_n(t)-1} 2 \sum_{l \in \N_0} \max\rbr{a_{m_n+l}-c,0}  = \sum_{n = n_0+1}^{+\ns}  2 k_{n}(t) \sum_{l \in \N_0} \max\rbr{a_{m_n+l}-c,0}  \nonumber \\
&  =  2 \sum_{n = n_0+1}^{+\ns}   2^{n} \zeta_t^2\sum_{l \in \N_0} \max\rbr{a_{m_n+l}-c,0}   \le \TTV {x^3}{[0,t]}c   \le \sum_{n\in \N_0} \sum_{k=0}^{k_n(t)} 2 \sum_{l \in \N_0} \max\rbr{a_{m_n+l}-c,0}  +1  \nonumber \\
 & =  \sum_{n\in \N_0}  2 ( k_{n}(t) +1) \sum_{l \in \N_0} \max\rbr{a_{m_n+l}-c,0} +1   \le 2\sum_{n = 0}^{n_0} 2^{n+1} \sum_{l \in \N_0} \max\rbr{a_{m_n+l}-c,0} \nonumber \\
& \quad +  2 \sum_{n = n_0+1}^{+\ns}  \sum_{l \in \N_0} (\max\rbr{a_{m_n+l}-c,0} +1)  +  2\sum_{n = n_0+1}^{+\ns}   2^{n} \zeta_t^2\sum_{l \in \N_0} \max\rbr{a_{m_n+l}-c,0}  +1 .  \label{calc32}
\end{align}
Similarly, as in Example 2, we prove that defining
\[
\varphi(c) := \frac{1}{1+ \TTV {x^3}{[0,1]}c}
\]
we get that
\[
\lim_{c \ra 0+} \varphi(c) \TTV {x^3}{[0,t]}c = \zeta^2_t.
\]
In this case however, by Theorem \ref{thm:meta1}, for any continuous $g:\R \ra \R$ we have
\begin{align}
\varphi(c)\int_{\R} \cross {y,c}{x^3}{[0,t]}  g(y) \dd y & \rightarrow \int_{(0, t]}g(x^3_{s-})\dd \zeta^2_s  = \int_{\mathcal C \cap (0, t]}g(x^3_{s-})\dd \zeta^2_s =  \int_{\mathcal C \cap (0, t]}g(x^3_{s})\dd \zeta^2_s \nonumber  \\
& \quad = \int_{\mathcal C \cap (0, t]}g\rbr{\zeta^2_s}\dd \zeta^2_s = \int_{(0, t]}g\rbr{\zeta^2_s}\dd \zeta^2_s = \int_{\sbr{0,{\zeta^2_t}}} g(s) \dd s . \label{localtim}
\end{align}

\textbf{Example 4.} If $p \in \N \setminus \cbr{1}$,  setting in Example 2 $$b_n = 2^{-\lfloor n/p \rfloor},$$ we get that the $p$th variation of $x^2$ along the sequence of Lebesgue partitions $\pi(\gamma) :=  \rbr{\pi\rbr{b_n, \gamma b_n}}_{n=1}^{+\ns}$, which we will denote by $\sbr{x^2}^{(p),\pi(\gamma)}_t$ and which is defined  as the limit $$\sbr{x^2}^{(p), \pi(\gamma)}_t := \lim_{n \ra +\ns}V^{b_n, \gamma b_n}_{p}\rbr{x^2, [0,t]},$$ where $V^{b_n, \gamma b_n}_{p}\rbr{x^2, [0,t]} := V^{b_n, \gamma b_n}_{\psi}\rbr{x^2, [0,t]}$ for $\psi(c) = c^p$ (recall $V^{c, \gamma c}_{\psi}$ defined in Subsection \ref{psivar}), exists and is equal to
\begin{equation} \label{pvarleb}
\sbr{x^2}^{(p),\pi(\gamma)}_t = 
\begin{cases}
\frac{2}{1- 2^{-(p-1)}} \zeta^2_t \text{ if } \gamma \in (0,1),\\
\frac{2\rbr{2^p-1}}{1- 2^{-(p-1)}} \zeta^2_t \text{ if } \gamma=0,
\end{cases}
\end{equation}
where $\zeta^2$ is the Cantor function defined in Example 2. As a consequence of \eqref{pvarleb} and \eqref{eeesstt} we obtain
\[
\lim_{n \ra +\ns} b_n^{p-1} \TTV {x^2}{[0,t]}{b_n} = \frac{2}{1- 2^{-(p-1)}} \zeta^2_t, \quad t \ge 0.
\]

To prove \eqref{pvarleb}, first notice that on each interval $ \bar{I}_{n,k}, n\in \N_0, k \in \cbr{0, 1, 2, \ldots, 2^n-1}$ and $N \in \N_0$, $N \ge n$, there are exactly $2 \lfloor \frac{b_n - \gamma b_N}{b_N} \rfloor$ times from the partition $\pi\rbr{b_N, \gamma b_N}$ and we have
$$2 \lfloor \frac{b_n - \gamma b_N}{b_N} \rfloor = 2 \lfloor \frac{2^{\lfloor \frac{N}{p} \rfloor}}{2^{\lfloor \frac{n}{p} \rfloor}} - \gamma \rfloor 
= \begin{cases}
2 \rbr{2^{\lfloor \frac{N}{p} \rfloor- \lfloor \frac{n}{p} \rfloor} } \text{ if } \gamma = 0 \text{ and } N \ge p \lfloor \frac{n}{p} \rfloor \\
2 \rbr{ 2^{\lfloor \frac{N}{p} \rfloor- \lfloor \frac{n}{p} \rfloor} - 1 } \text{ if } \gamma \in (0, 1) \text{ and } N \ge p \lfloor \frac{n}{p} \rfloor \\
0 \text{ otherwise}
\end{cases}.
$$
This, for $t \in I_{n_0,k}$, $n_0 \in \N$ and $\gamma \in (0,1)$ gives the following estimates 
\begin{equation}
 \label{mmm4p}
V^{b_N, \gamma b_N}_{p}\rbr{x^2, [0,t]} \ge 
b_N^p \sum_{n=p\lfloor n_0/p \rfloor +p}^{ p \lfloor N/p \rfloor + p -1} 2 k_n(t) \rbr{2^{\lfloor \frac{N}{p} \rfloor- \lfloor \frac{n}{p} \rfloor} -1}
\end{equation} 
and
\begin{equation} \label{mmm3p}
V^{b_N, \gamma b_N}_{p}\rbr{x^2, [0,t]}  \le b_N^p  \sum_{n=0}^{ p \lfloor n_0/p \rfloor + p -1} 2^{2n_0+2}  \rbr{2^{\lfloor \frac{N}{p} \rfloor- \lfloor \frac{n}{p} \rfloor}}  
+ b_N^p \sum_{n=p\lfloor n_0/p \rfloor +p}^{ p \lfloor N/p \rfloor + p -1} 2 k_n(t) \rbr{2^{\lfloor \frac{N}{p} \rfloor- \lfloor \frac{n}{p} \rfloor} -1}, 
\end{equation}
where $k_n(t)$ is defined in \eqref{kndef}. Using \eqref{calc23} and the fact that for $t \in I_{n_0,k}$  
\[
 \zeta_t^2 = \frac{1}{2^{n_0+1}} k_{n_0+1}(t) = \frac{1}{2^{n_0+2}} k_{n_0+2}(t) = \frac{1}{2^{n_0+3}} k_{n_0+3}(t) = \ldots
\]
we calculate 
\begin{align} 
& b_N^p \sum_{n=p\lfloor n_0/p \rfloor +p}^{ p \lfloor N/p \rfloor + p -1} 2 k_n(t) 2^{\lfloor \frac{N}{p} \rfloor- \lfloor \frac{n}{p} \rfloor}  = 2^{-p\lfloor \frac{N}{p} \rfloor} \sum_{n=p\lfloor n_0/p \rfloor +p}^{ p \lfloor N/p \rfloor + p -1} 2\cdot 2^{n} \zeta_t^2 \rbr{2^{\lfloor \frac{N}{p} \rfloor- \lfloor \frac{n}{p} \rfloor} } \nonumber \\
& = 2 \zeta_t^2 2^{-(p-1)\lfloor \frac{N}{p} \rfloor}  \sum_{n=p}^{ p \lfloor N/p \rfloor + p -1} 2^{n - \lfloor \frac{n}{p} \rfloor}  - A\rbr{n_0}=2 \zeta^2_t 2^{-(p-1)\lfloor \frac{N}{p} \rfloor} \rbr{1+ 2 + \ldots + 2^{p-1}} \sum_{m=1}^{\lfloor \frac{N}{p} \rfloor} 2^{p\cdot m - m} -  A\rbr{n_0} \nonumber \\
& =2 \zeta^2_t  2^{-(p-1)\lfloor \frac{N}{p} \rfloor} \rbr{2^p-1} 2^{p-1} \frac{2^{(p-1) \lfloor \frac{N}{p} \rfloor } - 1 }{2^{p-1}-1} -  A\rbr{n_0}  =2 \zeta^2_t  \frac{\rbr{2^p-1} 2^{p-1}}{2^{p-1}-1} \rbr{1 - 2^{-(p-1)\lfloor \frac{N}{p} \rfloor}  } -  A\rbr{n_0},
\label{mmm1p}
\end{align}
where 
$$A\rbr{n_0}= 2 \zeta^2_t 2^{-(p-1)\lfloor \frac{N}{p} \rfloor} \sum_{n=p}^{p\lfloor n_0/p \rfloor +p - 1} 2^{n - \lfloor \frac{n}{p} \rfloor}.
$$
We also have 
\begin{align}
& b_N^p  \sum_{n=p\lfloor n_0/2 \rfloor+p}^{ p \lfloor N/p \rfloor + p- 1} 2 k_n(t)   = 2^{-p \lfloor  \frac{N}{p} \rfloor}
 \sum_{n=p\lfloor n_0/2 \rfloor+p}^{ p \lfloor N/p \rfloor + p- 1}  2  \cdot 2^{n} \zeta_t^2 \nonumber  \\
& \quad =  2 \zeta^2_t  2^{-  p\lfloor  \frac{N}{p} \rfloor} \rbr{ 2^{p \lfloor \frac{N}{p} \rfloor +p}  - 2^{p\lfloor  \frac{n_0}{p} \rfloor + p}}  =  2 \zeta^2_t \rbr{ 2^{p}  - 2^{p\lfloor  \frac{n_0}{p} \rfloor + p -  p\lfloor  \frac{N}{p} \rfloor}}. \label{mmm2p}
\end{align}
Finally, using \eqref{mmm4p}-\eqref{mmm2p} and noticing that 
\[
 b_N^p  \sum_{n=0}^{ p \lfloor n_0/p \rfloor + p -1} 2^{2n_0+2}  \rbr{2^{\lfloor \frac{N}{p} \rfloor- \lfloor \frac{n}{p} \rfloor}} = 2^{-(p-1)\lfloor \frac{N}{p} \rfloor } \sum_{n=0}^{ p \lfloor n_0/p \rfloor + p -1} 2^{2n_0+2 - \lfloor \frac{n}{p} \rfloor }
\]
we obtain
\[
\lim_{N \ra +\ns} V^{b_N, \gamma b_N}_{p}\rbr{x^2, [0,t]} = 2 \zeta^2_t  \frac{\rbr{2^p-1} 2^{p-1}}{2^{p-1}-1} - 2 \zeta^2_t  2^{p}=    \frac{2}{1- 2^{-(p-1)}} \zeta^2_t.
\]
Similarly, we conclude that if $\gamma = 0$ then 
\[
\lim_{N \ra +\ns} V^{b_N, 0}_{p}\rbr{x^2, [0,t]} = \frac{2\rbr{2^p-1}}{1- 2^{-(p-1)}} \zeta^2_t.
\]
This gives \eqref{pvarleb}.

The same $p$-th variation equal the Cantor function (up to a multiplicative constant)  is obtained by Kim in \cite[Example 2.7]{Kim:2022} for $p \in 2\N$ who generalized Example 3.6 from Davis et al. \cite{Obloj_local:2015} for quadratic variation ($p=2$). Both examples were inspired by a remark on p. 194 in \cite{Bertoin:1987}, but the sequence of Lebesgue partitions used in \cite[Example 2.7]{Kim:2022} and \cite[Example 3.6]{Obloj_local:2015} is more complicated (different for different $I_n^k$, $n=0,1, \ldots$) than ours.

\section{Truncated variation and local times of self-similar processes with stationary increments} \label{Secthree}

Let us start with recalling that a process $X_t$, $t \in [0,+\ns)$, is called \emph{self-similar with index} $\beta >0$ if for any $A>0$
\begin{equation}
\rbr{ A^{-\beta}X_{At},\ t \in [0, +\ns) } =^{d}\rbr{X_{t},\ t \in [0, +\ns)},\label{eq:scalingProperty}
\end{equation}  
where $=^d$ denotes the equality of distributions. 

A stochastic process $X_t$, $t \in [0,+\ns)$,  \emph{has stationary increments} if for any $0\leq t_{1}<t_{2}\ldots<t_{n}$ and any $h>0$ we
have 
\begin{align*}
 & (X_{t_{2}}-X_{t_{1}}, X_{t_{3}}-X_{t_{2}}, \ldots, X_{t_{n}}-X_{t_{n-1}})\\
 & =^{d}(X_{t_{2}+h}-X_{t_{1}+h}, X_{t_{3}+h}-X_{t_{2}+h}, \ldots, X_{t_{n}+h}-X_{t_{n-1}+h}).
\end{align*}
\\{\bf Proof of Theorem \ref{thm:ssProcesses}:} By the definition of truncated variation \eqref{tvc_def} and (\ref{eq:scalingProperty}) for any $c>0$ we get 
\begin{align*}
&c^{1/\beta-1}\text{TV}^{c}(X,[0,t])=^{d}c^{1/\beta-1}\text{TV}^{c}(cX_{c^{-1/\beta} \cdot },[0,t])=c^{1/\beta}\text{TV}^{1}(X_{c^{-1/\beta} \cdot},[0,t])=c^{1/\beta}\text{TV}^{1}(X,[0,c^{-1/\beta}t]),
\end{align*}
where the last two equalities hold on the process level. We define the family
$\cbr{\eta(n,m)}_{m>n\geq0}$ of random variables 
\[
\eta(n,m):=\TTV X{[n,m]}1.
\]
The truncated variation depends only on the increments, hence 
$$\rbr{\eta(k,2k), \eta(2k,3k), \ldots}=^d \rbr{\eta(0, k), \eta(k, 2k), \ldots}, \quad k \in \N;$$
$$\rbr{\eta(k,k+1), \eta(k,k+2), \ldots}=^d \rbr{\eta(0, 1), \eta(0, 2), \ldots}, \quad k \in \N.$$
By superadditivity  \eqref{superadd} of the truncated variation we have 
\[
\eta(0,n)\geq\eta(0,m)+\eta(m,n).
\]
Recalling (\ref{eq:finiteExpectation}) we see that the double-indexed sequence $\rbr{\eta(n, m)}_{0\le n < m}$ satisfies assumptions of Kingman's Subadditive Ergodic Theorem in the version from \cite[Theorem 10.22]{Kallenberg:2001},
and we obtain 
\[
\lim_{n\conv+\infty}\frac{\eta(0,n)}{n}=Z  \text{ a.s.},
\]
where $Z\in(-\infty,+\infty]$ is a random variable measurable with respect to the invariant $\sigma$-field, which, by the ergodicity assumption, is trivial. Thus $Z$ is a.s. a constant, $Z = C\in [0, +\ns]$ a.s. Further, from \eqref{superaddd} we know
that 
\[
\eta(0,n)\leq\sum_{i=0}^{n-1}\eta(i,i+1)+n-1,
\]
hence we get $C<+\infty$. 

Thus we got the a.s. convergence of $c^{1/\beta}\text{TV}^{1}(X,[0,c^{-1/\beta}t])$ to a constant. This yields the convergence in probability of $c^{1/\beta-1}\text{TV}^{c}(X,[0,t])$ to the same constant since it has the same distribution as $c^{1/\beta}\text{TV}^{1}(X,[0,c^{-1/\beta}t])$.
\hfill  $\square$
\begin{cor}\label{cor:ssProcesses}
Under the assumptions of Theorem  \ref{thm:ssProcesses}, there exists a sequence of positive reals $\rbr{c_n}$, $n \in \N$, such that $c_n \ra 0+$ as $n \ra +\ns$, and a measurable set $\Omega_X \subset \Omega$ such that $\P \rbr{\Omega_X} = 1$ and for each $\omega \in \Omega_X$ and  each $t\in[0, +\ns)$
\[
c_n^{1/\beta-1}\TTVemph {X(\omega)}{[0, t]}{c_n} \ra C \cdot t \text{ as } n \ra +\ns,
\]
where $C$ is the same as in Remark  \ref{rem:constant}.
\end{cor}
\noindent{\bf Proof:} 
Let $\rbr{r_n}$, $n \in \N$, denote the sequence of all positive rational numbers. By Theorem \ref{thm:ssProcesses} there exists a sequence $\rbr{c_{1,n}}$ such that $c_{1,n} \ra 0+$ as $n \ra +\ns$ and a measurable set $\Omega_1 \subset \Omega$ such that $\P \rbr{\Omega_1} = 1$ and for each $\omega \in \Omega_1$
\[
c_{1,n}^{1/\beta-1}\TTV {X(\omega)}{[0, r_1]}{c_{1,n}} \ra C \cdot r_1.
\]  
Next, there exists a subsequence $\rbr{c_{2,n}}$ of $\rbr{c_{1,n}}$ and a measurable set $\Omega_2 \subset \Omega$ such that $\P \rbr{\Omega_2} = 1$ and for each $\omega \in \Omega_2$
\[
c_{2,n}^{1/\beta-1}\TTV {X(\omega)}{[0, r_2]}{c_{2,n}} \ra C \cdot r_2.
\]  
By the diagonal procedure, taking $c_n := c_{n,n}$ and $\Omega_X := \bigcap_n \Omega_n$ we obtain a sequence and a measurable set $\Omega_X \subset \Omega$ such that $\P \rbr{\Omega_X} = 1$ and for each $\omega \in \Omega_X$
\[
c_n^{1/\beta-1}\TTV {X(\omega)}{[0, t]}{c_n} \ra C \cdot t \text{ for each rational } t \in[0, +\ns).
\]
Now the convergence for all $t \in [0, +\ns)$ follows from the fact that $t \mapsto \TTV {x}{[0, t]}{c}$ is non-decreasing.
\hfill  $\square$

Using Theorem \ref{thm:meta1}, Remark \ref{subsequence}, 
and Corollary \ref{cor:ssProcesses} we will prove Corollary \ref{cor:ltssProcesses}.

\noindent{\bf Proof of Corollary \ref{cor:ltssProcesses}:} 
By Corollary \ref{cor:ssProcesses} there exists a sequence $\rbr{c_n}$, $n \in \N$, such that $c_n \ra 0+$ as $n \ra +\ns$ and a measurable set $\Omega_X \subset \Omega$ such that $\P \rbr{\Omega_X} = 1$ and for each $\omega \in \Omega_X$ and  each $t\in[0, +\ns)$
\begin{equation} \label{convvv}
c_n^{1/\beta-1}\TTV {X(\omega)}{[0, t]}{c_n} \ra C \cdot t \text{ as } n \ra +\ns,
\end{equation}
where $C$ is the same as in Theorem  \ref{thm:ssProcesses} and Remark  \ref{rem:constant}.
Let $g: \R \ra \R$ be a continuous function and $\omega \in \Omega_X$. Using \eqref{convvv}, Theorem \ref{thm:meta1} and Remark \ref{subsequence} we have that 
\begin{align} 
 c_n^{1/\beta-1} \int_{\R} \cross{z,c_n}{X(\omega)}{[0, t]} g(z)\dd z & \rightarrow C \int_{(0, t]}g(X_{u-}(\omega))\dd u  = C \int_{(0, t]}g(X_{u}(\omega))\dd u. \label{first_conv}
\end{align}
On the other hand, using the occupational times formula, we have
\begin{equation}  \label{occtfor}
\int_{(0, t]} g(X_{u}(\omega))\dd u = \int_{\R} g(y)  L_t^y(\omega) \dd y. 
\end{equation}
Comparing \eqref{first_conv} and \eqref{occtfor} we obtain that
\begin{align*} 
 c_n^{1/\beta-1} \int_{\R} \cross{z,c_n}{X(\omega)}{[0, t]} g(z)\dd z & \rightarrow C \int_{\R} g(y)  L_t^y(\omega) \dd y,
\end{align*}
which yields the claimed weak convergence.
\hfill  $\square$

By now we obtained almost sure weak convergence of the measures $$c_n^{1/\beta-1} \crossemph{y ,c_n}{X(\omega)}{[0, t]} \dd y$$ for some unidentified sequence $\rbr{c_n}_{n \in \N}$.
To formulate the next result let us introduce the notion of \emph{weak, locally uniform convergence in probability}.
\begin{dfn}
Let $(\Omega, {\cal F}, \P)$ be a probability space and let $\mu_t^c$,  $t \ge 0$, $c \ge 0$, be a family of finite,  random measures on $(\R, {\cal B}(\R))$. We will say that the random measures $\mu_t^c$ tend \emph{weakly locally uniformly in probability} as $c \ra 0+$ to the random measures $\mu_t^0$, $t \ge 0$, if for any continuous and bounded $g: \R \ra \R$ and any $T \ge0$ the supremum
\[
D^c_T(\omega) :=\sup_{t \in [0, T]} \left| \int_{\R} g \dd \mu_t^c(\omega) - \int_{\R} g \dd \mu_t^0(\omega) \right|
\]
is measurable relative to ${\cal F}$ and ${\cal B}(\R)$, thus is a random variable, and $D_T^c$ tend in probability $\P$ to $0$ as $c \ra 0+$.
\end{dfn}
Now we are ready to state the next result as follows. 
\begin{cor}\label{cor:luwssProcesses}
Assume that a real process $X_t$, $t \in [0, +\ns)$, satisfies assumptions of Theorem  \ref{thm:ssProcesses} and possesses a local time $L$ (in the sense of Definition \ref{localt_process}) relative to the Lebesgue measure $\dd y$ on $\R$ and the constant mapping $\Xi(\omega) \equiv \dd u$.  Then the measures
\[
c^{1/\beta-1} \crossemph{y ,c}{X(\omega)}{[0, t]} \dd y, \quad t\ge 0,
\]
tend weakly locally uniformly in probability as $c \ra 0+$ to the measures $C\cdot L_t^y \dd y$, where $C$ is the same as in Remark  \ref{rem:constant}.\end{cor}
\noindent{\bf Proof:} 
Let $g: \R \ra \R$ be a continuous, bounded function, let $T \ge 0$ and let us consider any sequence of positive reals $\rbr{c_n}$, $n \in \N$, converging to $0$. We need to prove that
\[
D_T^n := \sup_{t \in [0, T]} \left| c_n^{1/\beta-1}  \int_{\R} \cross{y ,c_n}{X(\omega)}{[0, t]} g(y )\dd y - \int_{\R} g(y) L_t^y(\omega) \dd y \right|, \quad n \in \N,
\]
are random variables which tend in probability to $0$ as $n \ra +\ns$. 

{The measurability of $D_T^n$ follows from the \cad property of $X$}. 

To prove the convergence in probability first notice that by the occupation times formula and the fact that there are only countable many jumps,
\begin{align*}
 \int_{\R} g(y) L_t^y(\omega) \dd y= \int_0^t g\rbr{X_{s}(\omega)} \dd s  =  \int_0^t g\rbr{X_{s-}(\omega)} \dd s .
\end{align*}
We will use the well known subsequence criterion, see \cite[Lemma 4.2]{Kallenberg:2001}: $D_T^n$ tends in probability to $0$ iff every subsequence of $\N$, $\rbr{d_n}_{n \in \N}$, has further subsequence $\rbr{d_{n_q}}$, $q \in \N$, such that $D_T^{d_{n_q}} $ tends a.s. to $0$.
Indeed, let $\rbr{d_n}_{n \in \N}$ be a subsequence of $\N$. Using the same diagonal procedure as in the proof of Corollary \ref{cor:ssProcesses} we obtain a further subsequence $\rbr{d_{n_q}}$, $q \in \N$, and a measurable set $\tilde{\Omega} \subset \Omega$ such that $\P \rbr{\tilde{\Omega}} = 1$ and for  each $t \in[0, +\ns)$ and  $\omega \in \tilde{\Omega}$
\[
c_{d_{n_q}}^{1/\beta-1}\TTV {X(\omega)}{[0, t]}{c_{d_{n_q}}} \ra C \cdot t.
\]
By Theorem \ref{thm:meta1}, Remark \ref{subsequence} and the fact that the linear function $[0, +\ns) \ni t \mapsto C \cdot t$ is continuous, we get that $D_T^{d_{n_q}} $ tends a.s. to $0$ as $q \ra +\ns$.

\hfill  $\square$

\section{Examples and related results} \label{secfour}

In this section, whenever we refer to  local time (without further elaboration) of a stochastic process, we mean local time in the sense of Definition \ref{localt_process} relative to the Lebesgue measure $\dd y$ on $\R$ and the constant mapping $\Xi(\omega) \equiv \dd u$, where $\dd u$  is the Lebesgue measure on $[0, +\ns)$.

\subsection{Fractional Brownian motions}

A standard example of a self-similar process with stationary increments where Corollaries \ref{cor:ltssProcesses}  and \ref{cor:luwssProcesses} may be applied is fractional Brownian motion. The fractional Brownian motion $B^H_t,$ $t \ge 0$, with the Hurst parameter $H \in \rbr{0,1}$, is a centred Gaussian process which is self-similar (with index $\beta = H$) and has stationary increments. More precisely,
\[
B^H_t - B^H_s \sim {\cal N}\rbr{0, \text{const}.|t-s|^{2H}} \text{ for } s, t\ge 0.
\]
Last, but not least -- $B^H$ possesses local time, see \cite{Berman:1973}.
\begin{fact} Let $B^H_t,$ $t \ge 0$, be a fractional Brownian motion with the Hurst parameter $H \in \rbr{0,1}$ and $L$ its local time.
There exists a sequence $\rbr{c_n}$, $n \in \N$, such that $c_n \ra 0+$ as $n \ra +\ns$ and a measurable set $\Omega_H \subset \Omega$ such that $\P \rbr{\Omega_H} = 1$ and for each $\omega \in \Omega_H$ and  each $t\in[0, +\ns)$
\begin{equation} \label{fbmweeaak}
c_n^{1/H-1} \crossemph{y ,c_n}{B^H(\omega)}{[0, t]} \dd y \ra^{weakly} C \cdot L_t^y(\omega) \dd y \text{ as } n \ra +\ns,
\end{equation}
where $C$ is the same as in Remark  \ref{rem:constant}. 

We also have that the measures
\[
c^{1/H-1} \crossemph{y ,c}{B^H(\omega)}{[0, t]} \dd y, \quad t\ge 0,
\]
tend weakly locally uniformly in probability as $c \ra 0+$ to the measures $C\cdot L_t^y \dd y$.
\end{fact}
\noindent{\bf Proof:} By Corollaries \ref{cor:ltssProcesses} and \ref{cor:luwssProcesses} it is sufficient to check the assumptions of Theorem  \ref{thm:ssProcesses} and the fact if $B^H$ possesses a local time $L$  relative to the Lebesgue measure $\dd y$ on $\R$ and the constant mapping $\Xi(\omega) \equiv \dd u$. 

Ergodicity of the sequence
$$\rbr{\TTV {B^H}{[0,k]}1, \TTV {B^H}{[k,2k]}1, \TTV {B^H}{[2k,3k]}1, \ldots }$$
follows from the mixing property of the increments of $B^H$, which may be verified using the autocorrelation function of the increments, see \cite{Ito:gauss1944}. More precisely, the quantities  $\TTV {B^H}{[lk,(l+1)k]}1$, $l \in \N_0$, may be approximated by functions of the increments $B^H_{lk+jk/N}-B^H_{lk+(j-1)k/N}$, $j=1,2,\ldots,N$, for sufficiently large (but fixed) $N \in \N$. But for $l, L \in \N_0$, $\cbr{B^H_{lk+ik/N}-B^H_{lk+(i-1)k/N}, i=1,2,\ldots,N}$ and $\cbr{B^H_{Lk+jk/N}-B^H_{Lk+(j-1)k/N}, j=1,2,\ldots,N}$ become independent when $l$ is fixed and $L \ra + \ns$ since it may be verified  by a direct calculation that then the correlation of $B^H_{lk+ik/N}-B^H_{lk+(i-1)k/N}$  
and $B^H_{Lk+jk/N}-B^H_{Lk+(j-1)k/N}$, $i,j = 1, 2,\ldots, N$, tends to $0$.

The finiteness of $\E \TTV {B^H}{[0, t]}{c}$ for any $c>0$ follows easily from the tail estimates of $\TTV {B^H}{[0, t]}{c}$ given in
\cite[Corollary 3]{BednorzLochowski:2013}. \hfill  $\square$

It is already mentioned in the Introduction, that a much stronger result for $H \in (0, 1/2)$ was proven in \cite{burdzy:2009, Toyomu:2023}. It states (see \cite[Theorem 1.3 or Theorem 1.4]{Toyomu:2023}) that for each $H \in (0, 1/2)$ there exists $\tilde{\Omega}_H \subset \Omega$ such that $\P \rbr{\tilde{\Omega}_H} = 1$ and for each $\omega \in \tilde{\Omega}_H$ and  each $t\in[0, +\ns)$
\[
\lim_{c \ra 0+} c^{1/H-1} \cross{y ,c}{B^H(\omega)}{[0, t]} = C \cdot L_t^y(\omega).
\]
The analogous result for $H \in (1/2, 1)$ is, as far as we know, not known, however, using \cite[Theorem 1.1]{Toyomu:2023} and Theorem \ref{thm:meta1} we can strengthen \eqref{fbmweeaak} and obtain the following fact.
\begin{fact} Let $B^H_t,$ $t \ge 0$, be a fractional Brownian motion with the Hurst parameter $H \in \rbr{0,1}$.
The convergence  \eqref{fbmweeaak} holds with probability $1$ for any sequence $\rbr{c_n}$, $n \in \N$, such that $c_n >0$, $n \in \N$, and 
$\lim_{n \ra + \ns} c_n = 0$. 
\end{fact}
\noindent{\bf Proof:}
Indeed, by \cite[Theorem 1.1]{Toyomu:2023} and relation \eqref{eeesstt} we get that
\[
\lim_{n \ra + \ns} n^{-\eta (1/H-1)}  \TTV {B^H}{[0, t]}{n^{-\eta}} = C \cdot t \text{ a.s.}
\]
where $\eta$ is some positive real. But this yields a similar convergence for any sequence $\rbr{c_n}$, $n \in \N$, such that $c_n >0$, $n \in \N$, and 
$\lim_{n \ra + \ns} c_n = 0$ because for $c_n <1$ we have the estimates $\lfloor 1/c_n^{1/\eta} \rfloor^{-\eta} \ge c_n$, $\lceil 1/c_n^{1/\eta} \rceil^{-\eta} \le c_n$ thus $\TTV {B^H}{[0, t]}{\lfloor 1/c_n^{1/\eta} \rfloor^{-\eta}} \le \TTV {B^H}{[0, t]}{c_n} $, $\TTV {B^H}{[0, t]}{\lceil 1/c_n^{1/\eta} \rceil^{-\eta} } \ge \TTV {B^H}{[0, t]}{c_n} $ and 
\begin{align*}
& \rbr{ \frac{\lceil 1/c_n^{1/\eta} \rceil}{\lfloor 1/c_n^{1/\eta} \rfloor}}^{-\eta (1/H-1)}   \lfloor 1/c_n^{1/\eta} \rfloor^{-\eta(1/H-1)} \TTV {B^H}{[0, t]}{\lfloor 1/c_n^{1/\eta} \rfloor^{-\eta}}\quad \le\quad c_n^{1/H-1} \TTV {B^H}{[0, t]}{c_n} \\
& \le\rbr{ \frac{\lfloor 1/c_n^{1/\eta} \rfloor}{\lceil 1/c_n^{1/\eta} \rceil}}^{-\eta (1/H-1)}  {\lceil 1/c_n^{1/\eta} \rceil}^{-\eta (1/H-1)} \TTV {B^H}{[0, t]}{{\lceil 1/c_n^{1/\eta} \rceil}^{-\eta}}.
\end{align*}
Since 
\[
\lim_{n \ra +\ns }{ \frac{\lfloor 1/c_n^{1/\eta} \rfloor}{\lceil 1/c_n^{1/\eta} \rceil}} = 1,
\]
we get 
that
\[
\lim_{c_n \ra 0+} c_n^{1/H-1} \TTV {B^H}{[0, t]}{c_n} = C \cdot t \text{ a.s.}.
\]
Now, applying Theorem \ref{thm:meta1} we have
\begin{equation*}
c_n^{1/H-1} \cross{y ,c_n}{B^H(\omega)}{[0, t]} \dd y \ra^{weakly} C \cdot L_t^y(\omega) \dd y \quad \text{ as } n \ra +\ns \text{ a.s.}
\end{equation*}
\hfill  $\square$

\subsection{Rosenblatt processes with the Hurst parameter $H \in (1/2,1)$}
Rosenblatt processes are neither Markov processes nor semimartingales nor Gaussian processes.  They live in the second Wiener chaos, in contrast to Gaussian processes (which belong to the first chaos). Similarly as fBms, they are parametrized by the Hurst parameter $H$. The parameter $H$ belongs to the interval $\rbr{1/2, 1}$. $R^H_t$, $t \in [0, +\ns)$,  with the Hurst parameter $H$ may be represented as 
\[
R^H_t = \text{const.} \int_{\R} \int_{\R} \rbr{\int_0^t \rbr{s - u_1}_+^{-\frac{2-H}{2}}\rbr{s - u_2}_+^{-\frac{2-H}{2} \dd s }}\dd B_{u_1} \dd B_{u_2},
\]
where $B$ is a standard Brownian motion. 
$R_H$ is self-similar with index $H$, has stationary increments, a.s. continuous paths and has the same autocorrelation function as a fBm with the same Hurst index, see \cite{Davis:2011}.
Rosenblatt processes possess local time, see \cite{Kerchev:2021}.
\begin{fact} Let $R^H_t,$ $t \ge 0$, be a Rosenblatt process with the Hurst parameter $H \in \rbr{1/2,1}$  and $L$ its local time.
There exists a sequence $\rbr{c_n}$, $n \in \N$, such that $c_n \ra 0+$ as $n \ra +\ns$ and a measurable set $\Omega_H \subset \Omega$ such that $\P \rbr{\Omega_H} = 1$ and for each $\omega \in \Omega_H$ and  each $t\in[0, +\ns)$
\begin{equation*} 
c_n^{1/H-1} \crossemph{y ,c_n}{R^H(\omega)}{[0, t]} \dd y \ra^{weakly} C \cdot L_t^y(\omega) \dd y \text{ as } n \ra +\ns,
\end{equation*}
where $C$ is the same as in Remark  \ref{rem:constant}. 

We also have that the measures
\[
c^{1/H-1} \crossemph{y ,c}{R^H(\omega)}{[0, t]} \dd y, \quad t\ge 0,
\]
tend weakly locally uniformly in probability as $c \ra 0+$ to the measures $C\cdot L_t^y \dd y$.
\end{fact}
\noindent{\bf Proof:} 
Rosenblatt processes have all properties assumed in Theorem \ref{thm:ssProcesses}. 

Ergodicity may be proven in a similar way as for fBms, using the fact that the finite-dimensional distributions of the increments become independent after large shifts, see also \cite{Kiska:2022}. 

Finiteness of 
 $\E \TTV {R^H}{[0, t]}{c}$ for any $c>0$ follows easily from
\cite[Corollary 2]{BednorzLochowski:2013} and tail estimates of the increments, see \cite[Proposition 4.2]{Kerchev:2021} or \cite{Major:2007}. 

Now we can use Corollaries \ref{cor:ltssProcesses} and \ref{cor:luwssProcesses}.
\hfill  $\square$

\subsection{Strictly $\alpha$-stable processes with $\alpha \in (1,2]$}

Other examples of self-similar processes satisfying assumptions of  Corollary \ref{cor:ltssProcesses} provide strictly $\alpha$-stable processes $X^{\alpha}$ with $\alpha \in (1,2]$. They are self-similar with index $\beta = 1/\alpha$, see \cite[Chapt. 15]{Kallenberg:2001}, and possess local time, see a general result of Hawkes in \cite[Chapt. V, Theorem 1]{Bertoin:1996uq}. 
\begin{fact} Let $X^{\alpha}_t,$ $t \ge 0$, be a strictly $\alpha$-stable process with $\alpha \in (1,2]$ and $L$ its local time.
There exists a sequence $\rbr{c_n}$, $n \in \N$, such that $c_n \ra 0+$ as $n \ra +\ns$ and a measurable set $\Omega_{\alpha} \subset \Omega$ such that $\P \rbr{\Omega_{\alpha}} = 1$ and for each $\omega \in \Omega_{\alpha}$ and  each $t\in[0, +\ns)$
\begin{equation*}
c_n^{\alpha-1} \crossemph{y ,c_n}{X^{\alpha}(\omega)}{[0, t]} \dd y \ra^{weakly} C \cdot L_t^y(\omega) \dd y \text{ as } n \ra +\ns,
\end{equation*}
where $C$ is the same as in Remark  \ref{rem:constant}. 

We also have that the measures
\[
c^{\alpha-1} \crossemph{y ,c}{X^{\alpha}(\omega)}{[0, t]} \dd y, \quad t\ge 0,
\]
tend weakly locally uniformly in probability as $c \ra 0+$ to the measures $C\cdot L_t^y \dd y$.
\end{fact}
\noindent{\bf Proof:} Ergodicity of the sequence
$$\rbr{\TTV {X^{\alpha}}{[0,k]}1, \TTV {X^{\alpha}}{[k,2k]}1, \TTV {X^{\alpha}}{[2k,3k]}1, \ldots }$$
 follows directly from the independence of the increments of $X^{\alpha}$.

The finiteness of $\E \TTV {X^{\alpha}}{[0, t]}{c}$ for any $c>0$ is proven in \cite{Cloud:2022} for the case  $\alpha \in (1,2)$ and in \cite{BednorzLochowski:2013} for the case  $\alpha =2$.

Now the fact follows from Corollaries \ref{cor:ltssProcesses} and \ref{cor:luwssProcesses}.
\hfill  $\square$
\begin{remark}
Notice that the normalization $c_n^{\alpha -1}$ of the numbers of crossings the intervals $\rbr{y-c_n/2, c_n/2}$ used to obtain the local time $L^y_t$ is the same as the normalization for the number of excursions of $X^{\alpha}$ which start from $y$ and hit the set $\left(-\ns, -c_n\right]\cup \left[c_n, +\ns\right)$, calculated in \cite[Example 6.3]{FristedTaylor:1983}. The proper normalization of the number of downcrossings from hitting $y+c/2$ to hitting $y-c/2$ needed to obtain the local time of more general L\'{e}vy processes is given in \cite{GetoorCrossLevy:1976}.

Compare the just obtained convergence of $$c_n^{\alpha-1} \crossemph{z,c_n}{X^{\alpha} }{[0, t]}$$  with a general result for c\`adl\`ag semimartingales obtained in \cite[Theorem 3.4]{LochOblPS:2021}, from which it follows that 
$$c_n \crossemph{z,c_n} {X^{\alpha} }{[0, t]} \xrightarrow[]{\text{a.s. in } \mathbb{L}^1} \mathcal L^z_t,$$ 
where $\mathcal L^z_t$ is a semimartingale local time, that is it satisfies the relation
\[
\int_0^t {\mathbf 1}_{\Gamma}(X_u^{\alpha}(\omega))  \dd [X^{\alpha}]_u^{cont} = \int_{\Gamma}\mathcal L_t^y(\omega) \dd y, \quad \Gamma \in {\cal B}(\R),
\]
where $[X^{\alpha}]^{cont}$ denotes the continuous part of the quadratic variation of $X^{\alpha}$. Since $[X^{\alpha}]^{cont} \equiv 0$ for $\alpha \in (1,2)$ and $[X^{2}]^{cont}_u \equiv \text{const.} u$, these results are not in contradiction.

\end{remark}

\subsection{Strictly $1$-stable process}
If $X$ is the strictly $1$-stable process (the Cauchy process) then the thesis of Theorem \ref{thm:ssProcesses} is not true since it is a self-similar process with index $\beta = 1$ and  $\int_{\R} \cross{y,c_n}{X^1(\omega)}{[0, t]} \dd y$ tends to $+\ns$ a.s. as $c_n \ra 0+$, since the trajectories of the Cauchy process have a.s. infinite total variation. 
Let us note that all the assumptions of Theorem \ref{thm:ssProcesses} are satisfied except the finiteness of $\E \TTV {X^{1}}{[0, t]}{c}$ for some $c>0$. Let us also note that if $X$ is strictly $1$-stable and symmetric then $X$ has no local time, see \cite[Chapt. V, Theorem 1]{Bertoin:1996uq}.

However, it is possible to prove (Piotr Mi\l o\'s, personal communication) that for almost all trajectories $x$ of $X^1$ Theorem \ref{thm:meta1} holds with the normalization $ \varphi(c) = {1}/{\ln (1/c)}$ for $c \approx 0$ and $\zeta_t = \text{const}. t$. The same normalization (in order to obtain the local time at the level $y$) needs to be applied for the number of excursions which start from $y$ and hit the set $\left(-\ns, -c\right]\cup \left[c, +\ns\right)$ of the \emph{asymmetric} Cauchy process, see \cite[Example 6.4]{FristedTaylor:1983}.

\subsection{Deterministic paths from Examples 1, 2 and 4 of Sect. \ref{sectwo}
}

\subsubsection{Deterministic path from Example 1 of Sect. \ref{sectwo}}

Naturally, the function $x^1$  has no local time in the sense of Definition \ref{def_local_det} relative to the measure $\dd \zeta^1$ and the Lebesgue measure $\dd y$ on $\R$ since $\dd \zeta^1$ is the Dirac delta,  $\dd \zeta^1= \delta_0$, and is not absolutely continuous with respect to the Lebesgue measure $\dd y$. 

\subsubsection{Deterministic path from Examples 2 and 4 of Sect. \ref{sectwo}}

Similarly, the function $x^2$  has no local time in the sense of Definition \ref{def_local_det} relative to the measure $\dd \zeta^2$ and the Lebesgue measure $\dd y$ on $\R$. If it was not the case, then for $g_{\varepsilon}: \R \ra [0,1]$ such that $g_{\varepsilon}$ is continuous, $g_{\varepsilon}(0) = 1$ and $g_{\varepsilon}(y) = 0$ for $ y \in \R \setminus (-\varepsilon, \varepsilon)$ and for some Lebesgue measurable function $L_1: \R \ra [0, + \ns)$ we would have 
\begin{align*}
1& =  \int_{\mathcal C} 1 \dd \zeta^2_s   = \int_{\mathcal C} g_{\varepsilon}\rbr{0} \dd \zeta^2_s  = \int_{\mathcal C} g_{\varepsilon}\rbr{x^2_s} \dd \zeta^2_s   = \int_{[0,1]} g_{\varepsilon}\rbr{x^2_s} \dd \zeta^2_s   = \int_{\R} g_{\varepsilon}(y) L^y_1 \dd y \ra_{{\varepsilon} \ra 0} 0.
\end{align*}

Let  $F$ belong to the space $C^p(\R, \R)$ (that is $F:\R \ra \R$ and has continuous $p$th derivative), $p \in 2 \N$, and let $\pi(\gamma)$ be the sequence of Lebesgue partition defined in Example 4. Then the corresponding change of variable formula 

\[
F\rbr{x^2_t} - F\rbr{x^2_0} = \int_0^t F'\rbr{x^2_s} \dd x^2_s + \frac{1}{p!} \int_0^t F^{(p)}\rbr{x^2_s} \dd \sbr{x^2}^{(p), \pi(\gamma)}_s
\]
of Cont and Perkowski holds (see \cite{ContPerkowski:2018}, \cite{Kim:2022}), but the corresponding local time $L_1: \R \ra [0, +\ns)$ of $x^2$ of order $p$ satisfying
\begin{equation} \label{last}
\int_0^1 {\mathbf 1}_{\Gamma}\rbr{x_s^{2}}\dd\sbr{x^2}^{(p), \pi(\gamma)}_s = \frac{1}{p!} \int_{\Gamma} L_1^y \dd y, \quad \Gamma \in {\cal B}(\R),
\end{equation}
does not exist. The proof is the same as in \cite[Example 2.7]{Kim:2022} and follows from the fact that the measure $\dd\sbr{x^2}^{(p), \pi(\gamma)} = \text{const}. \zeta^2$ does not charge the intervals $I_{n,k}$, $n \in \N_0$, $k=0,1,2,\ldots, 2^n-1$, defined in Example 2, thus for any $L_1$ satisfying \eqref{last} we would have $\int_{[b,+\ns)} L_1^y \dd y = 0$ for any $b >0$.

\subsection{The deterministic path $x^3$ from Example 3 of Sect. \ref{sectwo}
}

\label{subsect46}

In this case, from \eqref{localtim} it follows that the function $x^3$  has local time in the sense of Definition \ref{def_local_det} relative to the measure $\dd \zeta^2$ and the Lebesgue measure $\dd y$ on $\R$. This local time is equal
\[
L_t^y = {\mathbf 1}_{\sbr{0,\zeta_t^2}}(y).
\]
Additionally, let us notice that this example also shows that the weak convergence of the measures 
$
\varphi(c)\cross {y,c}{x^3}{[0,t]} \dd y 
$
to 
$L_t^y \dd y$,
stemming from Theorem \ref{thm:meta1} and the existence of the local time, can not be replaced by the weak convergence of $\varphi(c)\cross {y,c}{x^3}{[0,1]}$ to $L_1^y = {\mathbf 1}_{\sbr{0,1}}(y)$ in $\mathbb{L}^1 (\R, \dd y)$, that is we have the following fact.
\begin{fact}
The functions $\varphi(c)\crossemph {\cdot ,c}{x^3}{[0,1]}$ do not tend in $\mathbb{L}^1 (\R, \dd y)$ to the function $L_1 = {\mathbf 1}_{\sbr{0,1}}(\cdot)$ in $\mathbb{L}^1 (\R, \dd y)$ as $c \ra 0+$.
\end{fact}
\noindent{\bf Proof:}
If it was not the case, for any Borel measurable and essentially bounded function  $g: \R \ra \R$ we would have
\[
\lim_{c \ra 0+} \varphi(c) \int_{\R} \cross {y,c}{x^3}{[0,1]} g(y) \dd y =  \int_{\R} {\mathbf 1}_{\sbr{0,1}}(y) g(y) \dd y = \int_0^1 g(y) \dd y.
\]
From \eqref{calc300} it follows that 
$$B = \bigcup_{n \in \N_0} \bigcup_{k = 0}^{2^{n}-1} x^3 \rbr{I_{n,k}}$$
has the Lebesgue measure no greater than
\[
\sum_{n \in \N_0} \sum_{k = 0}^{2^{n}-1} a_{m_n} = \sum_{n \in \N_0} 2^n  a_{m_n} <1.
\]
Further, let us notice that for any $y \in [0,1]\setminus B$ and any $c>0$ we have $ \cross {y,c}{x^3}{[0,1]}  = 1$ (the function $x^3$ does not decrease below the value $y$ after attaining it). Thus, taking $g = {\mathbf 1}_{[0,1]\setminus B}$ we get 
\begin{align*}
& \lim_{c \ra 0+} \varphi(c) \int_{\R} \cross {y,c}{x^3}{[0,1]} g(y) \dd y = \lim_{c \ra 0+} \varphi(c) \int_{\R} {\mathbf 1}_{[0,1]\setminus B}(y)\dd y 
= \lim_{c \ra 0+} \varphi(c) \int_{[0,1]\setminus B} 1 \dd y = 0
\end{align*}
but on the other hand, if the claimed weak convergence in $\mathbb{L}^1 (\R, \dd y)$ was true, we would have
\begin{align*}
& \lim_{c \ra 0+} \varphi(c) \int_{\R} \cross {y,c}{x^3}{[0,1]} g(y) \dd y  = \int_{0}^1 {\mathbf 1}_{[0,1]\setminus B}(y)\dd y 
 \ge 1 - \sum_{n \in \N_0} 2^n a_{m_n} > 0.
\end{align*}
\hfill $\square$
\section*{Acknowledgments} The authors are very grateful to the anonymous reviewer of this paper, whose friendly, useful and detailed comments improved the presentation of the results of this paper. The authors are also grateful to Nicolas Perkowski who read a part of an earlier version of this paper, for his comments. The research of WMB and RM{\L} was supported by grant no. 2022/47/B/ST1/02114 \emph{Non-random equivalent characterizations of sample boundedness} of National Science Centre, Poland.


\begin{thebibliography}{0}

\bibitem[Ber73]{Berman:1973}
S.~M. Berman, \emph{Local nondeterminism and local times of gaussian process},
  Bull. Am. Math. Soc. \textbf{79} (1973), no.~2, 475--477.

\bibitem[Ber87]{Bertoin:1987}
J.~Bertoin, \emph{Temps locaux at int\'egration stochastique pour les processus
  de {D}irichlet.}, S\'eminaire de {P}robabilit\'es {XXI}, Lecture Notes in
  Math., vol. 1247, Springer, Berlin, 1987, pp.~191--205.

\bibitem[Ber96]{Bertoin:1996uq}
Jean Bertoin, \emph{L\'evy processes}, Cambridge Tracts in Mathematics, vol.
  121, Cambridge University Press, Cambridge, 1996. MR{1406564 (98e:60117)}

\bibitem[B{\L}15]{BednorzLochowski:2013}
W.~M. Bednorz and R.~M. {\L}ochowski, \emph{Integrability and concentration of
  sample paths' truncated variation of fractional {B}rownian motions,
  diffusions and {L}{\'e}vy processes}, Bernoulli \textbf{21} (2015), no.~1,
  437--464.

\bibitem[B{\L}M22]{Cloud:2022}
Witold~M. Bednorz, Rafa{\l}~M. {\L}ochowski, and Rafa{\l} Martynek, \emph{On
  optimal uniform approximation of {L{\'e}vy} processes on {Banach} spaces with
  finite variation processes}, ESAIM, Probab. Stat. \textbf{26} (2022),
  378--396 (English).

\bibitem[BY14]{BertoinYor:2014}
Jean Bertoin and Marc Yor, \emph{Local times for functions with finite
  variation: two versions of Stieltjes change-of-variables formula}, Bulletin
  of the London Mathematical Society \textbf{46} (2014), no.~3, 553 -- 560.

\bibitem[CLJPT81]{Chacon:1981}
R.~V. Chacon, Y.~Le~Jan, E.~Perkins, and S.~J. Taylor, \emph{Generalised arc
  length for Brownian motion and L{\'e}vy processes}, Zeitschrift f{\"u}r
  Wahrscheinlichkeitstheorie und Verwandte Gebiete \textbf{57} (1981), no.~2,
  197--211.

\bibitem[CP19]{ContPerkowski:2018}
R.~Cont and N.~Perkowski, \emph{Pathwise integration and change of variable
  formulas for continuous paths with arbitrary regularity}, Trans. Amer. Math.
  Soc. Ser. B \textbf{6} (2019), 161--186.

\bibitem[D{\L}MP23]{Toyomu:2023}
P.~Das, R.~M. {\L}ochowski, T.~Matsuda, and N.~Perkowski, \emph{Level crossings
  of fractional brownian motion}, Preprint arXiv:2308.08274 (2023).

\bibitem[DOS18]{Obloj_local:2015}
M.~Davis, J.~Ob{\l }\'{o}j, and P.~Siorpaes, \emph{Pathwise stochastic calculus
  with local times}, Ann. de l'IHP \textbf{54(1)} (2018), 1--21.

\bibitem[EK78]{ElKaroui:1978}
N.~El~Karoui, \emph{Sur les mont\'{e}es des semi-martingales}, Ast\'{e}risque
  \textbf{52--53 (Temps Locaux )} (1978), 63--88.

\bibitem[FT83]{FristedTaylor:1983}
Bert Fristed and S.~J. Taylor, \emph{Constructions of local time for a Markov
  process}, Z. Wahrscheinlichkeitstheorie verw. Gebiete \textbf{62} (1983),
  73--112.

\bibitem[Get76]{GetoorCrossLevy:1976}
R.~K. Getoor, \emph{Another limit theorem for local time}, Zeitschrift f{\"u}r
  Wahrscheinlichkeitstheorie und Verwandte Gebiete \textbf{34} (1976), no.~1,
  1--10.

\bibitem[GH80]{GemanHorowitz80}
Donald Geman and Joseph Horowitz, \emph{{Occupation densities}}, The Annals of
  Probability \textbf{8} (1980), no.~1, 1 -- 67.
  
\bibitem[HM{\L}Z24]{Hove:2024}
D.~Hove, F.~J. Nhlanga,  R.~M. {\L}ochowski, and P.~L. Zondi, \emph{Local times of deterministic paths with finite variation}, Preprint arXiv:2405.13174 (2024).


\bibitem[It{\^o}44]{Ito:gauss1944}
Kiyosi It{\^o}, \emph{{On the ergodicity of a certain stationary process}},
  Proceedings of the Imperial Academy \textbf{20} (1944), no.~2, 54 -- 55.

\bibitem[Kal01]{Kallenberg:2001}
O.~Kallenberg, \emph{Foundations of modern probability}, Probability and Its
  Applications, Springer, Heidelberg, Berlin, 2001.

\bibitem[Kis22]{Kiska:2022}
B.~Ki\v{s}ka, \emph{Variation of fractional processes}, Master Thesis, Charles
  University, {https://dspace.cuni.cz/handle/20.500.11956/171647} (2022).

\bibitem[Kim22]{Kim:2022}
Donghan Kim, \emph{Local times for continuous paths of arbitrary regularity},
  Journal of Theoretical Probability \textbf{35} (2022), no.~4, 2540--2568.

\bibitem[KNSV21]{Kerchev:2021}
George Kerchev, Ivan Nourdin, Eero Saksman, and Lauri Viitasaari, \emph{Local times and sample path properties of the Rosenblatt process}, Stochastic
  Processes and their Applications \textbf{131} (2021), 498--522.

\bibitem[L{\'e}vy40]{Levy:1940}
P.~L\'{e}vy, \emph{Le mouvement brownien plan}, Amer. J. Math. \textbf{62}
  (1940), 487--550.

\bibitem[Lem83]{Lemieux:1983}
M.~Lemieux, \emph{On the quadratic variation of semimartingales}, Master
  Thesis, The University of British Columbia,
  {https://circle.ubc.ca/handle/2429/23964} (1983).

\bibitem[{\L}17]{LochowskiColloquium:2017}
R.~M. {\L}ochowski, \emph{On a generalisation of the {B}anach {I}ndicatrix {T}heorem}, Colloq. Math.
\textbf{148(2)} (2017), 301--314.

\bibitem[{\L}G14]{LochowskiGhomrasni:2014}
R.~M. {\L}ochowski and R.~Ghomrasni, \emph{Integral and local limit theorems
  for level crossings of diffusions and the skorohod problem}, Electron. J.
  Probab. \textbf{19(10)} (2014), 1--33.

\bibitem[{\L}G15]{LochowskiGhomrasniMMAS:2015} \emph{The play operator, the truncated variation and the generalisation of the {J}ordan decomposition}, Math. Methods Appl. Sci.
  \textbf{38} (2015), no.~3, 403--419.

\bibitem[{\L}OPS21]{LochOblPS:2021}
Rafa{\l}~M. {\L}ochowski, Jan Ob{\l}\'oj, David~J. Pr\"omel, and Pietro
  Siorpaes, \emph{Local times and tanaka--meyer formulae for c{\`{a}}dl{\`{a}}g
  paths}, Electron. J. Probab. \textbf{26} (2021), 1--29.

\bibitem[Maj07]{Major:2007}
Peter Major, \emph{{On a Multivariate Version of Bernstein's Inequality}},
  Electronic Journal of Probability \textbf{12} (2007), no.~none, 966 -- 988.

\bibitem[Per81]{Perkins:1981}
Edwin Perkins, \emph{A global intrinsic characterization of Brownian local
  time}, The Annals of Probability \textbf{9} (1981), no.~5, 800--817.

\bibitem[PP15]{PerkowskiProemel_local:2015}
N.~Perkowski and D.~J. Pr\"{o}mel, \emph{Local times for typical price paths
  and pathwise {T}anaka formulas}, Electron. J. Probab. \textbf{20(46)} (2015),
  1--15.
 
\bibitem[BKR09]{burdzy:2009}
K. Burdzy , W. Kang and K. Ramanan,  \emph{The Skorokhod problem in a time-dependent interval}. {\em Stochastic Processes And Their Applications}. \textbf{119}, 428-452 (2009)


\bibitem[RY05]{RevuzYor:2005}
Daniel Revuz and Marc Yor, \emph{Continuous martingales and {B}rownian motion,
  3rd ed.}, Grundlehren der Mathematischen Wissenschaften, vol. 293,
  Springer-Verlag, Berlin, 2005. {MR1083357 (92d:60053)}

\bibitem[Taq11]{Davis:2011}
Murad Taqqu, \emph{The {R}osenblatt process in: Selected works of {M}urray
  {R}osenblatt ({R}ichard {A}. {D}avis, {K}eh-{S}hin {L}ii and {D}imitris {N}.
  {P}olitis, eds.)}, pp.~29--45, Springer New York, New York, NY, 2011.

\end{thebibliography}

\bibliographystyle{imsart-number}

\end{document}